# INVARIANCE PRINCIPLE FOR THE RANDOM CONDUCTANCE MODEL WITH UNBOUNDED CONDUCTANCES

By M. T. Barlow[1] and J.-D. Deuschel[2]

*University of British Columbia and Technische Universität Berlin*

We study a continuous time random walk $X$ in an environment of i.i.d. random conductances $\mu_e \in [1, \infty)$. We obtain heat kernel bounds and prove a quenched invariance principle for $X$. This holds even when $\mathbb{E}\mu_e = \infty$.

**1. Introduction.** We consider the Euclidean lattice $\mathbb{Z}^d$ with $d \geq 2$. Let $E_d$, the set of nonoriented nearest neighbour bonds, and, writing $e = \{x, y\} \in E_d$, let $(\mu_e, e \in E_d)$ be nonnegative r.v., defined on a probability space $(\Omega, \mathbb{P})$. Throughout this paper we will assume that $(\mu_e)$ is stationary and ergodic, and that its law is invariant under symmetries of $\mathbb{Z}^d$. We write $\mu_{xy} = \mu_{\{x,y\}} = \mu_{yx}$, and let $\mu_{xy} = 0$ if $x \not\sim y$. Set

$$\mu_x = \sum_y \mu_{xy}, \qquad P(x, y) = \frac{\mu_{xy}}{\mu_x}. \tag{1.1}$$

There are two natural continuous time random walks associated with $\mu$. Both jump according to the transitions $P(x, y)$. The first (the *constant speed random walk* or *CSRW*) waits at $x$ for an exponential time with mean 1 while the second (the *variable speed random walk* or *VSRW*) waits at $x$ for an exponential time with mean $1/\mu_x$. Write $\mathcal{L}_C$ and $\mathcal{L}_V$ for their generators, given by

$$\mathcal{L}_C f(x) = \mu_x^{-1} \sum_y \mu_{xy}(f(y) - f(x)), \tag{1.2}$$

$$\mathcal{L}_V f(x) = \sum_y \mu_{xy}(f(y) - f(x)). \tag{1.3}$$

Received October 2008.
[1]Supported in part by NSERC (Canada) and EPSRC (UK).
[2]Supported in part by DFG-Forschergruppe 718.

Set

$$\mathcal{E}(f,g) = \frac{1}{2} \sum_{x \in \mathbb{Z}^d} \sum_{y \in \mathbb{Z}^d} \mu_{xy}(f(x) - f(y))(g(x) - g(y)).$$

Let $\nu_x = 1$, $x \in \mathbb{Z}^d$. It is easy to check that if $f$, $g$ have finite support, then

(1.4) $$\mathcal{E}(f,g) = -\sum_x g(x) \sum_y \mu_{xy}(f(y) - f(x)),$$

and so

(1.5) $$\mathcal{E}(f,g) = -\langle \mathcal{L}_V f, g \rangle_\nu = -\langle \mathcal{L}_C f, g \rangle_\mu.$$

Thus the VSRW is the Markov process associated with the Dirichlet form $(\mathcal{E}, \mathcal{D}(\mathcal{E}))$ on $L^2(\nu)$ and has stationary measure $\nu$ while the CSRW is the Markov process associated with the Dirichlet form $(\mathcal{E}, \mathcal{D}(\mathcal{E}))$ on $L^2(\mu)$ and has stationary measure $\mu$.

Let $X = (X_t, t \geq 0, P_\omega^x, x \in \mathbb{Z}^d)$ be either the CSRW or the VSRW. Write $\mathcal{L}$ for its generator, $\theta$ for its invariant measure (so either $\theta = \nu$ or $\theta = \mu$) and let

(1.6) $$q_t^\omega(x,y) = \frac{P_\omega^x(X_t = y)}{\theta_y}$$

be the transition density of $X$ (or heat kernel associated with $\mathcal{L}$). This model, of a reversible (or symmetric) random walk in a random environment, is often called the *random conductance model* or RCM, particularly in the special case when $(\mu_e)$ are i.i.d. We are interested in the long-range behavior of $X$ and, in particular, in obtaining heat kernel bounds for $q_t^\omega(x,y)$ and a quenched or $\mathbb{P}$-a.s. invariance principle for $X$. When $\mathbb{E}\mu_e < \infty$, an averaged invariance principle is obtained in [17].

We begin by discussing the case when $(\mu_e)$ are i.i.d. If $\mu_e = 0$ then $X$ never jumps across $e$. So if $p_+ = \mathbb{P}(\mu_e > 0)$ is less than $p_c(d)$, the critical probability for bond percolation in $\mathbb{Z}^d$, then $X$ is $\mathbb{P} \times P_\omega^x$-a.s. confined to a finite set. Thus we restrict to the case $p_+ > p_c$. A number of different authors have studied this model under various restrictions on the support of $\mu_e$. If $\mu_e \in \{0,1\}$ then this problem reduces to that of a random walk on (supercritical) percolation clusters (see [1] for heat kernel bounds, and [10, 29, 34] for quenched invariance principles). More generally it is useful to consider the following special cases:

Case 0. $c^{-1} \leq \mu_e \leq c$ for some $c \geq 1$;

Case 1. $0 \leq \mu_e \leq 1$;

Case 2. $1 \leq \mu_e < \infty$.

For case 0, heat kernel bounds follow from the results in [18, 19], and a quenched invariance principle is proved in [34]. Case 1 is treated in [11, 12, 30]. (The papers [11, 12] consider a discrete time random walk.) These papers prove an invariance principle, with a strictly positive diffusion constant $\sigma^2$. Further, [11] shows that Gaussian heat kernel bounds do not hold in general in this case.

In this paper we will look at case 2. There is not a great difference between the CSRW and VSRW in case 1, but in case 2, and in particular when $\mathbb{E}\mu_e = \infty$, the VSRW and CRSW do have different behaviors. Also, while the discrete time random walk with jump probabilities $P(x, y)$ given by (1.1) behaves in a similar fashion to the CSRW, there is no simple discrete time analogue of the VSRW in case 2. We remark that our result for the CSRW also gives an invariance principle for the discrete time random walk with jump probabilities $P(x, y)$.

Let

(1.7) $$X_t^{(\varepsilon)} = \varepsilon X_{t/\varepsilon^2}, \qquad t \geq 0.$$

Our first main result is the following quenched functional central limit theorem (QFCLT):

THEOREM 1.1. *Let $d \geq 2$. Suppose that $(\mu_e)$ are i.i.d., and $\mu_e \geq 1$ $\mathbb{P}$-a.s.*
(a) *Let $X$ be the VSRW. Then $\mathbb{P}$-a.s. $X^{(\varepsilon)}$ converges (under $P_\omega^0$) in law to a Brownian motion on $\mathbb{R}^d$ with covariance matrix $\sigma_V^2 I$ where $\sigma_V > 0$ is nonrandom.*
(b) *Let $X$ be the CSRW. Then $\mathbb{P}$-a.s. $X^{(\varepsilon)}$ converges (under $P_\omega^0$) in law to a Brownian motion on $\mathbb{R}^d$ with covariance matrix $\sigma_C^2 I$ where*

$$\sigma_C^2 = \begin{cases} \sigma_V^2/(2d\mathbb{E}\mu_e), & \text{if } \mathbb{E}\mu_e < \infty, \\ 0, & \text{if } \mathbb{E}\mu_e = \infty. \end{cases}$$

We also have heat kernel bounds for the VSRW.

THEOREM 1.2. *Let $d \geq 2$. Suppose that $(\mu_e)$ are i.i.d. and $\mu_e \geq 1$ $\mathbb{P}$-a.s. Let $q_t^\omega(x, y)$ be the heat kernel for the VSRW. Let $\eta \in (0, 1)$. There exist r.v. $U_x, x \in \mathbb{Z}^d$, such that*

(1.8) $$\mathbb{P}(U_x(\omega) \geq n) \leq c_1 \exp(-c_2 n^\eta)$$

*and constants $c_i$ (depending on $d$ and the distribution of $\mu_e$) such that the following hold.*
(a) *For all $x$, $y$, $t$*

(1.9) $$q_t^\omega(x, y) \leq c_3 t^{-d/2}.$$

(b) *If $|x - y| \vee t^{1/2} \geq U_x$, then*

(1.10) $$q_t^\omega(x, y) \leq c_3 t^{-d/2} e^{-c_4 |x-y|^2/t} \qquad \text{when } t \geq |x - y|,$$

(1.11) $$q_t^\omega(x, y) \leq c_3 \exp(-c_4 |x - y|(1 \vee \log(|x - y|/t)))$$
$$\text{when } t \leq |x - y|.$$

(c) Let $x, y \in \mathbb{Z}^d$, $t > 0$. Then

$$(1.12) \qquad q_t^\omega(x,y) \geq c_5 t^{-d/2} e^{-c_6 |x-y|^2/t} \qquad \text{if } t \geq U_x^2 \vee |x-y|^{1+\eta}.$$

(d) Let $x, y \in \mathbb{Z}^d$ and $t \geq c_7 \vee |x-y|^{1+\eta}$. Then

$$(1.13) \qquad c_8 t^{-d/2} e^{-c_9 |x-y|^2/t} \leq \mathbb{E} q_t^\omega(x,y) \leq c_{10} t^{-d/2} e^{-c_{11}|x-y|^2/t}.$$

Using Theorems 1.1 and 1.2 we can obtain a parabolic Harnack inequality (PHI) for $q_t^\omega$ (for the VSRW) by the same arguments as in [3] (see Theorem 4.7). Since the CSRW is a time change of the VSRW, harmonic functions and Green's functions are the same for both processes. The PHI for the VSRW implies an elliptic Harnack inequality (see Corollary 4.8) which therefore holds for both CSRW and VSRW. Combining the invariance principle and the PHI, we obtain, using the methods of [3], a local limit theorem for the VSRW (see Theorem 5.14).

When $d \geq 3$ the calculations in Section 6 of [3] then give bounds on the Green's function $g_\omega(x,y)$ defined by

$$(1.14) \qquad g_\omega(x,y) = \int_0^\infty q_t^\omega(x,y)\,dt.$$

THEOREM 1.3.  *Let $d \geq 3$, and suppose that $(\mu_e)$ are i.i.d. and $\mu_e \geq 1$ $\mathbb{P}$-a.s.*

(a) *There exist constants $c_1, \ldots, c_4$ and r.v. $U_x$, $x \in \mathbb{Z}^d$ such that*

$$(1.15) \qquad \mathbb{P}(U_x \geq n) \leq c_1 \exp(-c_2 n^{1/3})$$

*and*

$$(1.16) \qquad \frac{c_3}{|x-y|^{d-2}} \leq g_\omega(x,y) \leq \frac{c_4}{|x-y|^{d-2}} \qquad \text{if } |x-y| \geq U_x \wedge U_y.$$

(b) *Let $C = \Gamma(\frac{d}{2} - 1)/(2\pi^{d/2} \sigma_V^2)$. For any $\varepsilon > 0$ there exists $M = M(\varepsilon, \omega)$ with $\mathbb{P}(M < \infty) = 1$ such that*

$$(1.17) \qquad \frac{(1-\varepsilon)C}{|x|^{d-2}} \leq g_\omega(0,x) \leq \frac{(1+\varepsilon)C}{|x|^{d-2}} \qquad \text{for } |x| > M(\omega).$$

(c) *We have, $\mathbb{P}$-a.s.,*

$$(1.18) \qquad \lim_{|x|\to\infty} |x|^{2-d} g_\omega(0,x) = \lim_{|x|\to\infty} |x|^{2-d} \mathbb{E} g_\omega(0,x) = C.$$

REMARK 1.4.  (b) and (c) in Theorem 1.3 use the QFCLT, which in turn uses the ergodic theorem. As we do not have any rate of convergence in the QFCLT this means that [unlike the r.v. $U_x$ in (a)] we have no control on the tail of the r.v. $M$ in (b).

The main difference between the RCM and the percolation case is the possibility of traps. Suppose $e = \{x, y\}$ is a bond with $\mu_e = K \gg 1$, and that all the other bonds $e'$ adjacent to $x$ and $y$ have $\mu_{e'} \simeq 1$. Then $P(x, y) \simeq 1 - c/K$, so $X$ will jump between $x$ and $y$ $O(K)$ times before leaving $\{x, y\}$, and thus the CSRW will be trapped for a time of order $K$ in $\{x, y\}$. However, for the VSRW each jump takes only $O(K^{-1})$, so the total time spent in $\{x, y\}$ is only $O(1)$. A similar effect will arise from finite clusters of bonds of high conductivity.

The presence of traps of this kind is why we have, when $\mathbb{E}\mu_e = \infty$, that the diffusion constant $\sigma_C^2$ for the CSRW is zero. In this case it is natural to ask if a different scaling will give a nontrivial limit. There is a connection here with "aging" (see [6–8]), and in [2] it is proved that if the tail distribution $\mathbb{P}(\mu_e > t) \sim t^{-\alpha}$, then $\varepsilon^\alpha X_{t/\varepsilon^2}$ converges to the "fractional kinetic" motion with parameter $\alpha$ (see [8]).

While we have written this paper for the case of i.i.d. conductances $\mu_e$, our arguments do not require the full strength of this. In case 0, when the conductances are bounded (both above and below), then uniform upper and lower Gaussian heat kernel estimates, as in Theorem 1.2, are well known (see [18]). It follows (see Remark 6.3) that Theorem 1.1. holds for any stationary ergodic environment. On the other hand, in the unbounded case 2, there exist stationary ergodic environments such that the VSRW can explode in finite time; for an example, see Remark 6.6 below.

For the Gaussian bounds in Theorem 1.2 we need to control the sizes of the clusters of high conductivity which is done by comparison between the graph metric $d(x, y)$ and a new metric $\widetilde{d}$ (introduced by Davies in [15, 16]) which is adapted to the structure of the random walk and satisfies $\widetilde{d}(x, y) \leq \mu_{xy}^{-1/2}$ when $x \sim y$. This new metric is constructed by a first passage percolation procedure, and in this paper we have used first passage percolation arguments to compare the two metrics (see Lemma 4.2). These arguments use estimates from [25] which in turn use the fact that $\mu_e$ are i.i.d. However, we could also have used a direct argument as in [12], Lemma 3.1 or [30], Lemma 5.3. Once we have the Gaussian bounds (with sufficiently good control on the tails of the r.v. $U_x$ in Theorem 1.2), the quenched invariance principle follows with no further hypotheses on $\{\mu_e, e \in E_d\}$ other than that it is stationary, symmetric and ergodic. Theorem 6.1 summarizes the general situation.

The structure of this paper is as follows. In Sections 2 and 3 we study a deterministic graph $\Gamma = (G, E)$ with edge weights $\mu_{xy}$. Under certain conditions (which are $\mathbb{P}$-a.s. satisfied by the VSRW on the i.i.d. RCM) we obtain heat kernel bounds in this setup. Our approach uses similar methods to those used in [1] for percolation clusters. However, in [1] the Carne–Varopoulos "long-range" bounds played an essential role at various points.

These bounds do not hold for the VSRW, and instead we use more general upper bounds obtained by Davies [15, 16], which are in terms of the metric $\widetilde{d}(x, y)$. The same metric $\widetilde{d}$ is also needed to control $P_\omega^x(\widetilde{d}(X_t, x) \geq \lambda t^{1/2})$ which is the key step in obtaining general Gaussian upper bounds. Similar bounds, in the context of weighted Laplacians on manifolds, are obtained by Grigor'yan [24]; here the metric $\widetilde{d}$ is the Riemannian metric. Very recently, and independently, Mourrat [31] has obtained upper bounds for the VSRW which in certain cases improve on those in Theorem 1.2.

Once one has upper bounds, lower bounds follow by the same arguments as in [1], Section 5 (see Section 3). In Section 4 we then prove Theorem 1.2.

Section 5 proves the invariance principle. We begin with the VSRW. The basic technique in the proof (as in many previous papers such as [10, 12, 17, 27–29]) is to associate with $X_t$ a process $Z_t$ on $\Omega = [1, \infty]^{\mathbb{E}^d}$ which is the environment seen from the random walk. More precisely, for each $x \in \mathbb{Z}^d$, let $T_x : \Omega \to \Omega$ be given by

$$T_x(\omega)(z, w) = \omega(z + x, w + x).$$

Assuming that $X_0 = 0$ we define

$$(1.19) \qquad Z_t(\omega) = T_{X_t(\omega)}(\omega).$$

One seeks to use the process $Z$ to construct the "corrector," that is, a map $\chi : \Omega \times \mathbb{Z}^d \to \mathbb{R}^d$ such that

$$(1.20) \qquad M_t(\omega) = X_t(\omega) - \chi(\omega, X_t(\omega)), \qquad t \geq 0,$$

is a $P_\omega^0$-martingale. Once one has constructed the corrector, showing the invariance principle for the rescaled martingale $\varepsilon M_{t/\varepsilon^2}$ is standard, and using results from [12, 34], the heat kernel estimates in Theorem 1.2, together with the sublinear growth of $\chi$, imply that

$$(1.21) \qquad \lim_{\varepsilon \to 0} \varepsilon \chi(\omega, X_{t/\varepsilon^2}) = 0 \qquad \text{in } P_\omega^0\text{-probability.}$$

However, the standard construction of the corrector is based on $L^2(\mathbb{P})$ calculus, which requires finiteness of the first moment of the conductance (see [17], page 816). In our case we wish to handle the case when $\mathbb{E}\mu_e = \infty$, and so we need an alternative approach. (We remark that if $\mathbb{E}\mu_e = \infty$ then it is not easy to find suitable function spaces on $\Omega$ which give a core for the Dirichlet form associated with $Z$.)

Our solution relies on discretization. We define $\widehat{X}_n = X_n$, $n \in \mathbb{Z}_+$, and consider the process

$$(1.22) \qquad \widehat{X}_t^{(\varepsilon)} = \varepsilon \widehat{X}_{\lfloor t/\varepsilon^2 \rfloor}.$$

We can control $\sup_{t \leq T} |X_t^{(\varepsilon)} - \widehat{X}_t^{(\varepsilon)}|$ (see Lemma 4.12), so an invariance principle for $X^{(\varepsilon)}$ will follow from one for $\widehat{X}^{(\varepsilon)}$. The process $\widehat{X}$ does not have bounded jumps—in fact it jumps anywhere in $\mathbb{Z}^d$ with positive probability. However, the long-range bounds on $q_t^\omega(x,y)$ in (1.11) give good control on these jumps, and, in particular, the bounds on $q_1^\omega(x,y)$ imply that

$$\mathbb{E}E_\omega^0|X_1|^2 < \infty, \tag{1.23}$$

which is the key $L^2$ condition on $\widehat{X}$ for the construction of the corrector. As we will see in Section 5, looking at the discrete time process does actually introduce some simplifications in the construction of the corrector $\chi$. In the end (see Remark 5.15) it will turn out that the "discrete time" corrector $\chi$ also satisfies (1.20).

Finally, a short Section 6 makes some remarks on more general environments, and gives an example (a one-sided spanning tree) where the process $X$ fails to be conservative, and so the invariance principle fails.

**2. Transition density upper bounds on a fixed graph.** Let $\Gamma = (G, E)$ be an infinite (deterministic) graph, $\mu_e, e \in E$, be bond conductances and $\nu$ be a measure on $G$. We make the following assumptions on $(G, E)$, $\mu$ and $\nu$.

ASSUMPTION 2.1. (1) $\Gamma$ is connected.
(2) The vertex degree is uniformly bounded,

$$|\{y : y \sim x\}| \leq C_D \qquad \text{for all } x \in G. \tag{2.1}$$

(3) $\mu_e > 0$ for all $e \in E$.
(4) There exists $C_M \geq 1$ such that

$$C_M^{-1} \leq \nu_x = \nu(\{x\}) \leq C_M \qquad \text{for all } x \in G. \tag{2.2}$$

The results of this section do not explicitly require a strictly positive lower bound on $\mu_e$; however, a later assumption [see Assumption 2.6(2)] will impose some control on the edges $e$ with $\mu_e$ small.

We write $\mu_{xy}$ for $\mu_{\{x,y\}}$, and set $\mu_{xy} = 0$ if $x \not\sim y$. Let $d(x,y)$ be the usual graph distance on $G$, and write

$$B(x,r) = \{y : d(x,y) < r\}. \tag{2.3}$$

Let $C_A < \infty$. We now construct, by a first passage percolation procedure, a second metric $\widetilde{d}$ on $G$ satisfying

$$(C_A^{-2} \vee \mu_{yy'})|\widetilde{d}(x,y) - \widetilde{d}(x,y')|^2 \leq 1 \qquad \text{for every } x \in G, y \sim y'. \tag{2.4}$$

We write $\widetilde{B}(x,r) = \{y : \widetilde{d}(x,y) < r\}$ for balls in the metric $\widetilde{d}$. (In this paper we can take $C_A = 1$, but for possible future extensions we treat the general case.) To construct $\widetilde{d}$ we assign waiting times

$$t(e) = C_A \wedge \mu_e^{-1/2}, \qquad e \in E, \tag{2.5}$$

and then take $\widetilde{d}(x,y)$ to be the shortest journey time between $x$ and $y$. More precisely,

$$\widetilde{d}(x,y) = \inf\left\{\sum_{i=1}^{n} t(e_i)\right\}, \tag{2.6}$$

where the infimum is taken over paths $(e_1, \ldots, e_n)$ from $x$ to $y$. Since we do not have a strictly positive uniform lower bound on $t(e)$, in general there may not be a minimizing path. However, such paths will a.s. exist when $t(e_i)$ are i.i.d. positive random variables.

LEMMA 2.2. *The metric $\widetilde{d}$ constructed above satisfies (2.4).*

PROOF. Let $e = \{y, z\}$; then $\widetilde{d}(x,z) \leq \widetilde{d}(x,y) + t(e)$, and using (2.5) gives (2.4). □

Recall that

$$\mathcal{E}(f,g) = \frac{1}{2}\sum_x \sum_{y \sim x} \mu_{xy}(f(y) - f(x))(g(y) - g(x)).$$

Let $\mu_x = \sum_{y \sim x} \mu_{xy}$, and extend $\mu$ to a measure on $G$. Then $\mathcal{E}(f,g)$ is defined for $f, g \in L^2(G, \mu)$.

Let $X = (X_t, t \in [0, \infty), P^x, x \in G)$ be the continuous time Markov chain on $G$ with generator

$$\mathcal{L}f(x) = \nu_x^{-1} \sum_y \mu_{xy}(f(y) - f(x)).$$

At this point we cannot exclude the possibility that $X$ may explode, and we write $\zeta$ for the lifetime of $X$. Let

$$\|f\|_{\mathcal{E}_1}^2 = \mathcal{E}(f,f) + \|f\|_{L^2(\nu)}^2$$

and $\mathcal{F}$ be the closure of the set of functions on $G$ with finite support with respect to $\|f\|_{\mathcal{E}_1}$. Then $X$ is the Markov process associated with the Dirichlet form $(\mathcal{E}, \mathcal{F})$ on $L^2(G, \nu)$ (see [22]). Let $q_t(x,y)$ be the transition density (heat kernel) of $X$ with respect to $\nu$:

$$q_t(x,y) = \frac{P^x(X_t = y)}{\nu_y}.$$

We begin by using the results of Davies [15, 16] to obtain long-range bounds on $q_t$. By Proposition 5 of [15], we have

$$q_t(x,y) \leq (\nu_x \nu_y)^{-1/2} \inf_{\psi \in L^\infty(G)} \exp(\psi(y) - \psi(x) + \Lambda(\psi)t), \tag{2.7}$$

where $\Lambda(\psi) = \sup_x b(\psi, x)$, and

$$(2.8) \qquad b(\psi, x) = \frac{1}{2\nu_x} \sum_{y \sim x} \mu_{xy}(e^{\psi(x)-\psi(y)} + e^{\psi(y)-\psi(x)} - 2).$$

THEOREM 2.3. *Assume $(G, E)$ and $\mu$ satisfy Assumption 2.1. There exist constants $c_1, \ldots, c_4$ (depending only on $C_A, C_D, C_M$) such that the following hold.*

(a) *If $x, y \in G$ and $\widetilde{D} = \widetilde{d}(x, y) \leq c_1 t$, then*

$$(2.9) \qquad q_t(x, y) \leq c_2 \exp(-c_3 \widetilde{D}^2/t).$$

(b) *If $x, y \in G$ and $\widetilde{D} = \widetilde{d}(x, y) \geq c_1 t$, then*

$$(2.10) \qquad q_t(x, y) \leq c_2 \exp(-c_4 \widetilde{D}(1 \vee \log(\widetilde{D}/t))).$$

PROOF. Fix $x_0, y_0 \in G$, let $t > 0$, and write $\widetilde{D} = \widetilde{d}(x_0, y_0)$. Let $\lambda > 0$ and set

$$\psi_\lambda(x) = -\lambda(\widetilde{D} \wedge \widetilde{d}(x_0, x)), \qquad b(\lambda) = \sup_x b(\psi_\lambda, x).$$

Let $x \in G$, $y \sim x$ and write $\widetilde{\mu}_{xy} = C_A^{-2} \vee \mu_{xy}$,

$$J_{xy} = \mu_{xy}(e^{\psi_\lambda(x) - \psi_\lambda(y)} + e^{\psi_\lambda(y) - \psi_\lambda(x)} - 2).$$

Then as $|\psi_\lambda(x) - \psi_\lambda(y)| \leq \lambda \widetilde{\mu}_{xy}^{-1/2}$, and $\cosh(x)$ is increasing on $[0, \infty)$,

$$(2.11) \qquad J_{xy} \leq 2\mu_{xy}(\cosh(\lambda \widetilde{\mu}_{xy}^{-1/2}) - 1) \leq 2\widetilde{\mu}_{xy}(\cosh(\lambda \widetilde{\mu}_{xy}^{-1/2}) - 1).$$

Using the power series for cosh we have that the right-hand side of (2.11) is decreasing in $\widetilde{\mu}_{xy}$, so

$$J_{xy} \leq C_A^{-2}(e^{C_A \lambda} + e^{-C_A \lambda} - 2).$$

Hence

$$b(\psi_\lambda, x) \leq \tfrac{1}{2} C_M C_D C_A^{-2}(e^{C_A \lambda} + e^{-C_A \lambda} - 2).$$

Let $f(x) = e^x + e^{-x} - 2$; then $b(\lambda) \leq c_7 f(C_A \lambda)$. Thus by (2.7), and writing $y = C_A \lambda$,

$$(2.12) \qquad \begin{aligned} q_t(x_0, y_0) &\leq C_M \inf_\lambda \exp(-\lambda \widetilde{D} + c_7 t f(C_A \lambda)) \\ &\leq C_M \exp\left(\frac{\widetilde{D}}{C_A}\left(\inf_{y > 0}\left(-y + \frac{C_A c_7 t}{\widetilde{D}} f(y)\right)\right)\right). \end{aligned}$$

So if
$$F(s) = \inf_{y>0}(-y + (2s)^{-1}(e^y + e^{-y} - 2)),$$

then

(2.13) $$q_t(x_0, y_0) \leq C_M \exp\left(\frac{\widetilde{D}}{C_A} F\left(\frac{\widetilde{D}}{2C_A c_7 t}\right)\right)$$

and it remains to bound $F$.

We have (see [15], page 70) that
$$F(s) = s^{-1}((1 + s^s)^{1/2} - 1) - \log(s + (1 + s^2)^{1/2})$$

and also that $F(s) \leq -s/2(1 - s^2/10)$ for $s > 0$. Hence, if $s \leq 3$, then $F(s) \leq -s/20$ while if $s \geq e$, then
$$F(s) \leq 1 - \log(2s) = -\log(2s/e).$$

Substituting in (2.13) completes the proof. □

REMARK 2.4. Note that if $\widetilde{D} = ct$ then both (2.9) and (2.10) give a bound of the form $ce^{-t/c}$.

Since $\mu_e$ are not bounded above, the process $X$ may explode. The following condition is enough to exclude this.

LEMMA 2.5. *Suppose there exists $x \in G$ and $\theta > 0$ such that*

(2.14) $$\sum_{y \in G} \exp(-\theta \widetilde{d}(x, y)) \nu_y < \infty.$$

*Then $X$ is conservative.*

PROOF. Let $\zeta$ be the lifetime of $X$. Then as $\Gamma$ is connected it is easy to see that either $P^y(\zeta = \infty) = 1$ for all $y \in G$, or else $P^y(\zeta < t) > 0$ for all $y \in G$, $t > 0$.

For $n \geq C_A^{-2}$ let $\mu_{xy}^{(n)} = n \wedge \mu_{xy}$, $X^{(n)}$ be the process associated with the conductances $\mu^{(n)}$, and $q_t^{(n)}(x, y)$ be the transition density of $X^{(n)}$ with respect to $\nu$. We have $q_t(x, y) = \lim_{n \to \infty} q_t^{(n)}(x, y)$. Note that each $X^{(n)}$ is conservative, and that the bounds in Theorem 2.3 hold (with the same constants $c_i$) for each $q^{(n)}$. With (2.2) the condition (2.14) implies that $\widetilde{B}(x, R)$ is finite for each $R > 0$.

Let $t > 0$. With constants $c_i$ as in Theorem 2.3, choose $r$ large enough so that $r > c_1 t$ and $c_4(1 \vee \log(r/t)) \geq \theta$. So, if $R \geq r$, using (2.10),

(2.15) $$\sum_{y \in \widetilde{B}(x,R)^c} q_t^{(n)}(x, y) \nu_y \leq \sum_{y \in \widetilde{B}(x,R)^c} c_2 \exp(-\theta \widetilde{d}(x, y)) \nu_y < \infty.$$

Let $\varepsilon > 0$; then we can take $R$ large enough so that the right-hand side of (2.15) is less than $\varepsilon$. Thus, as $X^{(n)}$ is conservative, for all $n$,

$$\sum_{y \in \widetilde{B}(x,R)} q_t^{(n)}(x,y)\nu_y > 1 - \varepsilon.$$

So,

$$P^x(\zeta > t) \geq \sum_{y \in \widetilde{B}(x,R)} q_t(x,y)\nu_y = \lim_{n \to \infty} \sum_{y \in \widetilde{B}(x,R)} q_t^{(n)}(x,y)\nu_y \geq 1 - \varepsilon.$$

Therefore $P^x(\zeta > t) = 1$ for all $t$, proving that $X$ is conservative. $\square$

We now make further assumptions on the graph $\Gamma$ and the conductances $\mu$. As we will see in Section 4, it is easy to check these for the random conductance model on $\mathbb{Z}^d$ when $\mu_e \geq 1$.

ASSUMPTION 2.6. (1) There exists $d \geq 1$ and $C_V < \infty$ such that

(2.16) $\qquad \nu(B(x,r)) \leq C_V r^d \qquad$ for all $x \in G, r \geq 1$.

(2) There exists a constant $C_N$ such that the following Nash inequality holds. If $f \in L^1(G,\nu) \cap L^2(G,\nu)$, then, writing $\|f\|_p$ for norms in $L^p(G,\nu)$,

(2.17) $\qquad \mathcal{E}(f,f) \geq C_N \|f\|_2^{2+4/d} \|f\|_1^{-4/d}.$

REMARK 2.7. Note that (2) above does place some restrictions on the edges $e$ with $\mu_e$ small. For example, taking $x \in G$ and $f = 1_x$, (2.17) gives

$$\sum_y \mu_{xy} \geq C_N \nu_x^{1-2/d}.$$

By [14] we have:

LEMMA 2.8. *Suppose (2.17) holds. Then*

(2.18) $\qquad q_t(x,y) \leq ct^{-d/2}, \qquad x,y \in G, t > 0.$

To obtain better control of $q_t(x,y)$ when $d(x,y)$ is large we need to compare the metrics $d$ and $\widetilde{d}$. Note first that $\widetilde{d}(x,y) \leq C_A \wedge \mu_{xy}^{-1/2} \leq C_A$ when $x \sim y$, so

(2.19) $\qquad \widetilde{d}(x,y) \leq C_A d(x,y), \qquad x,y \in G.$

DEFINITION 2.9. Let $\lambda \geq 1$, $\eta \in (0,1)$. Let $x \in G$, $r \in [1,\infty)$. We say $(x,r)$ is $\lambda$-*good* if $\widetilde{B}(x, n/\lambda) \subset B(x,n)$ for all $n \geq r$, $n \in \mathbb{N}$. We say $(x, R_0)$ is $\lambda$-*very good* if for all $R \geq R_0$, $(y,r)$ is $\lambda$-good for all $y \in B(x, R)$, $r \geq R^\eta$, $r \in \mathbb{N}$. Note that if $(x, R_0)$ is $\lambda$-*very good* then $(x, R_1)$ is $\lambda$-*very good* for all $R_1 \geq R_0$. For $x \in G$ let $V_x = V_x(\lambda)$ be the smallest integer such that $(x, V_x)$ is $\lambda$-very good.

*Note*. Unlike the definitions in [1], the event that $(x,r)$ is $\lambda$-good depends on the structure of $\Gamma$ "at infinity."

LEMMA 2.10. *Suppose $(x, R)$ is $\lambda$-good.*
(a) *If $d(x, y) \geq R$,*

$$\lambda^{-1} d(x,y) \leq \widetilde{d}(x,y) \leq C_A d(x,y). \tag{2.20}$$

(b) *If $R' \geq (2R) \vee 2(1 + C_A \lambda) d(x, x')$, then $(x', R')$ is $2\lambda$-good.*

PROOF. (a) The upper bound is given in (2.19). For the lower bound, let $s = d(x, y) \geq R$. Then $y \notin B(x, s)$, so $y \notin \widetilde{B}(x, s/\lambda)$ and thus $\widetilde{d}(x, y) \geq s/\lambda$.
(b) Let $\widetilde{r} = \widetilde{d}(x, x')$, $r = d(x, x')$, and $s \geq R'$. Then as $s/2 \geq R$,

$$\widetilde{B}(x', s/2\lambda) \subset \widetilde{B}(x, \widetilde{r} + s/2\lambda) \subset B(x, \lambda \widetilde{r} + s/2).$$

So, using (2.20), $\widetilde{B}(x', s/2\lambda) \subset B(x', (1 + \lambda C_A) r + s/2) \subset B(x', s)$. □

LEMMA 2.11. *Let $x \in G$, $\theta > 0$, $r \geq 1$. If $(x,r)$ is $\lambda$-good, then*

$$\sum_{y \in B(x,r)^c} \exp(-\theta \widetilde{d}(x,y)) \nu_y \leq \begin{cases} c(\lambda) r^d e^{-cr\theta}, & \text{if } r\theta \geq 1, \\ c(\lambda) \theta^{-d}, & \text{if } r\theta < 1. \end{cases} \tag{2.21}$$

*In particular, $X$ is conservative.*

PROOF. Write $I$ for the left-hand side of (2.21), and $D_n = B(x, 2^n r) - B(x, 2^{n-1} r)$. Then

$$I \leq \sum_{n=1}^{\infty} \sum_{y \in D_n} \exp(-\theta d(x,y)/\lambda) \nu_y \leq \sum_{n=1}^{\infty} C_V (2^n r)^d \exp(-2^{n-1} r\theta/\lambda). \tag{2.22}$$

If $\alpha > 0$, $d \geq 1$ then there exists $c_1 = c_1(d)$ such that

$$\sum_{n=1}^{\infty} 2^{nd} e^{-\alpha 2^n} \leq \begin{cases} c_1 e^{-\alpha}, & \text{if } \alpha \geq 1, \\ c_1 \alpha^{-d}, & \text{if } \alpha < 1, \end{cases}$$

and using these bounds in (2.22) completes the proof. □

LEMMA 2.12. *Let $x \in G$ and suppose that $(x,r)$ is $\lambda$-good. If $t \in (0,1)$, then*

$$E^x d(x, X_t)^p \leq c(\lambda) r^{p+d}. \tag{2.23}$$

PROOF. Using the bound (2.10) a similar calculation to that in Lemma 2.11 gives

$$E^x d(x, X_t)^p \leq r^p + \sum_{n=1}^{\infty} C_V (2^n r)^{d+p} e^{-c\widetilde{d}(x,y)}$$

$$\leq r^p + c r^{d+p} \sum_{n=1}^{\infty} 2^{n(d+p)} \exp(-c' 2^n r / \lambda) \leq c r^{d+p}. \qquad \square$$

We now follow the arguments in [1], Section 3 (the "Bass–Nash method") to obtain Gaussian upper bounds on $q_t(x,y)$. As in Lemma 2.5 for $1 \leq n \leq \infty$, let $\mu_{xy}^{(n)} = \mu_{xy} \wedge n$, $X^{(n)}$ be the associated VSRW, and $q^{(n)}(x,y)$ be the transition density of $X^{(n)}$. Let $x_0 \in G$, and set

$$M_n(t) = M_n(x_0, t) = E^{x_0} \widetilde{d}(x_0, X_t^{(n)}) = \sum_y \widetilde{d}(x_0, y) q_t^{(n)}(x_0, y) \nu_y, \tag{2.24}$$

$$Q_n(t) = Q_n(x_0, t) = -\sum_y q_t^{(n)}(x_0, y) \log q_t^{(n)}(x_0, y) \nu_y. \tag{2.25}$$

The following three inequalities lead, by Nash's argument [32], to upper bounds on $M_n(t)$ which are uniform in $n$. [We remark that we only need the approximations $X^{(n)}$ to justify an interchange of sums in part (c) of the following lemma.]

LEMMA 2.13. *Let $x_0 \in G$ and $r \geq 1$. Suppose $(x_0, r)$ is $\lambda$-good, and $1 \leq n < \infty$. There exist constants $c_i$, independent of $n$, such that the following hold.*

(a) *We have, for $t > 0$,*

$$Q_n(x_0, t) \geq c_1 + \tfrac{1}{2} d \log t. \tag{2.26}$$

(b)

$$M_n(x_0, t) \geq c_2 e^{Q_n(x_0, t)/d} \qquad \textit{provided either } M_n(x_0, t) \geq r \textit{ or } t \geq c_3 r^2. \tag{2.27}$$

(c) *For $t > 0$,*

$$Q_n'(t) \geq c_4 M_n'(t)^2. \tag{2.28}$$

PROOF. We write $Q_n(t)$, $M_n(t)$ for $Q_n(x_0,t)$, $M_n(x_0,t)$. (a) is immediate from (2.18) and the fact that since $X^{(n)}$ is conservative, $\sum_y q_t^{(n)}(x_0,y)\nu_y = 1$.

(b) The proof is similar to those in [1, 4, 32]. First note that (2.16) and Lemma 2.8 give that

(2.29) $$M_n(t) \geq r \quad \text{provided } t \geq c_3 r^2.$$

By Lemma 2.11, provided $ar \leq 1$,

$$\sum_{y \in G} e^{-a\widetilde{d}(x_0,y)}\nu_y \leq \sum_{y \in B(x_0,r)} e^{-a\widetilde{d}(x_0,y)}\nu_y + \sum_{y \notin B(x_0,r)} e^{-a\widetilde{d}(x_0,y)}\nu_y$$

$$\leq cr^d + ca^{-d} \leq ca^{-d}.$$

Now $u(\log u + \theta) \geq -e^{-1-\theta}$ for $u > 0$. So, setting $\theta = a\widetilde{d}(x_0,y) + b$, where $a \leq 1/r$,

$$-Q_n(t) + aM_n(t) + b = \sum_{y \in G} q_t^{(n)}(x_0,y)(\log q_t^{(n)}(x_0,y) + a\widetilde{d}(x_0,y) + b)\nu_y$$

$$\geq -\sum_{y \in G} e^{-1-a\widetilde{d}(x_0,y)-b}\nu_y$$

$$\geq -e^{-1-b}\sum_{y \in G} e^{-a\widetilde{d}(x_0,y)}\nu_y \geq -c_5 e^{-b}a^{-d}.$$

Setting $a = 1/M_n(t)$ and $e^b = M_n(t)^d = a^{-d}$, we obtain

$$-Q_n(t) + 1 + d\log M_n(t) \geq -c_6$$

and rearranging gives (b).

(c) Set $f_t(x) = q_t^{(n)}(x_0,x)$, and let $b_t(x,y) = f_t(x) + f_t(y)$. We have, using (2.4),

$$M_n'(t) = \sum_y \widetilde{d}(x_0,y)\frac{\partial f_t(y)}{\partial t}\nu_y = \sum_y \widetilde{d}(x_0,y)\mathcal{L}^{(n)}f_t(y)\nu_y$$

$$= -\frac{1}{2}\sum_x \sum_y \mu_{xy}(\widetilde{d}(x_0,y) - \widetilde{d}(x_0,x))(f_t(y) - f_t(x));$$

the final interchange of sums can be justified using (2.9) and (2.10) and the fact that $\mu^{(n)}$ is uniformly bounded. So using (2.4),

$$M_n'(t) \leq \frac{1}{2}\sum_x \sum_{y \sim x}(\mu_{xy}^{1/2}|\widetilde{d}(x_0,y) - \widetilde{d}(x_0,x)|b_t(x,y)^{1/2})\left(\mu_{xy}^{1/2}\frac{|f_t(y) - f_t(x)|}{b_t(x,y)^{1/2}}\right)$$

$$\leq \frac{1}{2}\sum_x \sum_{y \sim x}(b_t(x,y)^{1/2})\left(\mu_{xy}^{1/2}\frac{|f_t(y) - f_t(x)|}{b_t(x,y)^{1/2}}\right)$$

$$\leq \frac{1}{2}\left(\sum_x \sum_{y \sim x} b_t(x,y)\right)^{1/2} \left(\sum_x \sum_y \mu_{xy} \frac{(f_t(y) - f_t(x))^2}{b_t(x,y)}\right)^{1/2}.$$

Now

(2.30) $$\sum_x \sum_{y \sim x} b_t(x,y) = 2 \sum_x \sum_{y \sim x} f_t(x) \leq 2 C_D C_M \sum_x f_t(x) \nu_x = 2 C_D C_M.$$

So,

$$M'_n(t)^2 \leq c \sum_x \sum_y \mu_{xy} \frac{(f_t(y) - f_t(x))^2}{f_t(x) + f_t(y)}.$$

Since we have, for $u, v > 0$,

$$\frac{(u-v)^2}{u+v} \leq (u-v)(\log u - \log v),$$

we deduce

$$M'_n(t)^2 \leq c \sum_x \sum_y \mu_{xy} (f_t(y) - f_t(x))(\log f_t(y) - \log f_t(x)).$$

Thus

$$Q'_n(t) = -\sum_y (1 + \log f_t(y)) \mathcal{L}^{(n)} f_t(y) \nu_y$$

(2.31) $$= \frac{1}{2} \sum_x \sum_y \mu_{xy} (\log f_t(y) - \log f_t(x))(f_t(y) - f_t(x))$$

$$\geq c M'_n(t)^2;$$

where again the interchange of sums uses (2.9) and (2.10) and the fact that $\mu^{(n)}$ is uniformly bounded. □

PROPOSITION 2.14. *Let $x_0 \in G$, $r \geq 1$ and $(x_0, r)$ be $\lambda$-good. Then*

(2.32) $$c_1 t^{1/2} \leq M_\infty(x_0, t) \leq c_2 t^{1/2} \qquad \text{for } t \geq c_3 r^2.$$

PROOF. Note that the lower bound is immediate from Lemma 2.8. For the upper bound let first $n < \infty$. Set $R_n(t) = d^{-1}(Q_n(t) - c_1 - \frac{1}{2} d \log t)$, so that $R_n(t) \geq 0$ by (2.26). Then

$$M_n(t) = \int_0^t M'_n(s) \, ds \leq c \int_0^t Q'_n(s)^{1/2} \, ds \leq c \int_0^t \left(R'_n(s) + \frac{1}{2s}\right)^{1/2} ds.$$

Using the inequality $(a+b)^{1/2} \leq b^{1/2} + a/(2b)^{1/2}$ and integrating by parts we obtain

$$M_n(t) \leq ct^{1/2} + c\int_0^t s^{1/2} R'_n(s)\, ds$$

$$\leq ct^{1/2} + c(1 + R_n(t)t^{1/2}) \leq ct^{1/2}(1 + R_n(t)).$$

By (2.27) we also have $M_n(t) \geq t^{1/2} e^{R_n(t)}$ for $t > cr^2$. Thus $R_n(t)$ is bounded for $t > cr^2$ and this implies that

(2.33) $$c_1 t^{1/2} \leq M_n(x_0, t) \leq c_2 t^{1/2} \qquad \text{for } t \geq c_3 r^2.$$

Since the constants in Lemma 2.13 are independent of $n$, the constants $c_i$ in (2.33) are also independent of $n$. Since $M_\infty(t) \leq \liminf_{n\to\infty} M_n(t)$, (2.32) then follows. $\square$

LEMMA 2.15. *Let $x \in G$, $r \geq 1$ and $(x,r)$ be $\lambda$-good. Then*

(2.34) $$c_1 t^{1/2} \leq E^x d(x, X_t) \leq c_2 t^{1/2} \qquad \text{for } t \geq c_3 r^2.$$

PROOF. Since $d(x, X_t) \geq C_A^{-1} \widetilde{d}(x, X_t)$, the first inequality is clear. Let $c_3$ be as in Proposition 2.14, and $t = c_3 R_1^2$, so $R_1 \geq r$. Then if $A = B(x, R_1)^c$, using (2.20),

$$E^x d(x, X_t) \leq \lambda R_1 + E^x(d(x, X_t); X_t \in A)$$

$$\leq \lambda R_1 + E^x(\lambda \widetilde{d}(x, X_t); X_t \in A)$$

$$\leq \lambda R_1 + \lambda c_4 t^{1/2} \leq c_5 t^{1/2}. \qquad \square$$

The next few results follow quite closely along the lines of [1]. Let

$$\tau(x, r) = \inf\{t : d(x, X_t) \geq r\}.$$

LEMMA 2.16. *There exist constants $c_1$, $c_2$, $c_3$ such that if $R \geq c_1$ and*

(2.35) $$(y, c_2 R) \text{ is } \lambda\text{-good for all } y \in B(x, R),$$

*then if $t_0 = R^2/(2c_3)$*

(2.36) $$P^x(\tau(x, R/2) \leq t_0) \leq \tfrac{1}{2}$$

*and hence for $t \geq 0$,*

(2.37) $$P^x(\tau(x, R) \leq t) \leq \frac{1}{2} + \frac{c_3 t}{R^2}.$$

PROOF. Write $\tau = \tau(x, R/2)$, and $c'_i$ for the constants $c_i$ in Lemma 2.15. Let $c_3 = 64(c'_2)^2$. Choose $c_2$ so that $r = c_2 R$ satisfies $c'_3 r^2 = t_0$, and $c_1$ so that $r \geq 1$. Then as (2.35) holds we can use Lemma 2.15 to bound $E^y d(y, X_s)$ for $s \geq t_0$, $y \in B(x, R)$. So,

$$c'_2 (2t_0)^{1/2} \geq E^x d(x, X_{2t_0}) \geq E^x (d(x, X_{t_0 \wedge \tau}) - d(X_{t_0 \wedge \tau}, X_{2t_0}))$$
$$\geq E^x 1_{(\tau \leq t_0)} d(x, X_\tau) - \sup_{y \in B(x,R)} \sup_{0 \leq s \leq t_0} E^y d(y, X_{2t_0 - s})$$
$$\geq P^x(\tau \leq t_0) R/2 - c'_2 (2t_0)^{1/2},$$

and rearranging we obtain (2.36).

Inequality (2.37) now follows easily; if $t \leq t_0$, $P^x(\tau(x, R) \leq t) \leq P^x(\tau(x, R/2) \leq t_0) \leq \frac{1}{2}$ while if $t > t_0$, then the right-hand side of (2.37) is greater than 1. □

To obtain the Gaussian upper bound on $q_t(x, y)$ we need to prove that the process $X$ does not move too rapidly across a ball $B(x, R)$. We choose $r \ll R$, and use the fact that if $X$ moves across $B(x, R)$ in the time interval $[0, t]$, then it has to move across many smaller balls of side $r$ in the same period; the estimate (2.37) is enough to bound the probability of this. Our argument uses the following easily proved estimate:

LEMMA 2.17 (See [5], Lemma 1.1). *Let $\xi_1, \xi_2, \ldots, \xi_n$, $V$ be nonnegative r.v. such that $V \geq \sum_1^n \xi_i$. Suppose that for some $p \in (0, 1)$, $a > 0$,*

(2.38) $$P(\xi_i \leq t | \sigma(\xi_1, \ldots, \xi_{i-1})) \leq p + at, \qquad t > 0.$$

*Then*

(2.39) $$\log P(V \leq t) \leq 2 \left( \frac{ant}{p} \right)^{1/2} - n \log(1/p).$$

PROPOSITION 2.18. *There exists constants $c_1, \ldots, c_4$ such that if $x \in G$, $R \geq c_1$, $t \geq c_1 R$ and*

(2.40) $$(z, c_2 t/R) \text{ is } \lambda\text{-good for all } z \in B(x, R),$$

*then*

(2.41) $$P^x(\tau(x, R) < t) \leq c_3 e^{-c_4 R^2 / t}.$$

PROOF. Let $1 \leq m < R/2$, and set $r = R/2m$. Define stopping times

$$S_0 = 0, \qquad S_i = \inf\{t \geq S_{i-1} : d(X_{S_{i-1}}, X_t) \geq r\}, \qquad i \geq 1.$$

Set $\tau_i = S_i - S_{i-1}$, and write $\mathcal{F}_t = \sigma(X_s, s \leq t)$ for the filtration of $X$. As $d(X_{S_i}, X_{S_{i+1}}) \leq r + 1 < 2r$, we have $S_m \leq \tau(x, R)$ and $X_{S_i} \in B(x, R)$ for $0 \leq i \leq m - 1$.

Suppose for the moment that $m$ is such that we can apply Lemma 2.16 to control each $\tau_i$. Then

$$(2.42) \qquad P^x(\tau_i \leq u | \mathcal{F}_{S_{i-1}}) \leq \frac{1}{2} + \frac{c_5 u}{r^2}, \qquad u > 0, 1 \leq i \leq m,$$

so writing $p = \frac{1}{2}$, $a = c_5/r^2$ and using (2.39), we obtain

$$\log P^x(\tau(x,r) < t) \leq \log P^x\left(\sum_1^m \tau_i < t\right)$$

$$(2.43) \qquad \leq 2(amt/p)^{1/2} - m\log p^{-1}$$

$$\leq -c_6 m\left(2 - \left(\frac{c_7 tm}{R^2}\right)^{1/2}\right)$$

$$= -c_6 m(2 - (m/\kappa)^{1/2}),$$

where $\kappa = R^2/(c_7 t)$. If $\kappa$ is such that we can choose $m \in \mathbb{N}$ with $\kappa \leq m < 2\kappa$, and so that (2.42) holds, then (2.43) implies (2.41).

We can choose $c_1$ so that $\kappa < R/2 - 1$. If $\kappa \leq 1$ then, adjusting the constant $c_3$ appropriately, (2.41) is immediate. If $1 < \kappa < R/2 - 1$ then let $m = \lfloor \kappa \rfloor + 1 \leq 2\kappa$. Then $\frac{1}{4}c_6(t/R) \leq r \leq \frac{1}{2}c_6(t/R)$, and so choosing $c_2$ suitably, (2.40) and Lemma 2.16 imply (2.42). $\square$

Recall that $V_x$ is the smallest $R$ such that $(x, R)$ is $\lambda$-very good.

THEOREM 2.19. *Let $x, y \in G$, and write $D = d(x,y)$. Suppose that either $D \geq c_1 \vee V_x$ or $t \geq D^2$.*

(a) *If $c_2 D \leq t$, then*

$$(2.44) \qquad q_t(x,y) \leq c_3 t^{-d/2} e^{-c_4 d(x,y)^2/t}.$$

(b) *If $c_2 D \geq t$, then*

$$(2.45) \qquad q_t(x,y) \leq c_3 \exp(-c_4 D(1 \vee \log(D/t))).$$

PROOF. We need to consider various cases. First, if $t \geq D^2$, then (2.44) follows from (2.18). So we can suppose $D \geq V_x$. Note that by (2.20) $\widetilde{d}(x,y) \geq \lambda^{-1} D$.

If $t \leq c_2 D$, then (2.45) follows from (2.10). If $c_2 D \leq t \leq c_6 D^2/\log D$, then Theorem 2.3 gives

$$q_t(x,y) \leq c_8 \exp(-2c_7 D^2/t).$$

Choosing $c_6$ small enough we have $\exp(-c_7 D^2/t) \leq t^{-d/2}$ and (2.44) follows.

It remains to consider the case $c_6 D^2 / \log D \le t \le D^2$. Let $A_x = \{z : d(x,z) \le d(y,z)\}$, $A_y = G - A_x$, $t' = t/2$, $D' = D/2$. Note that $B(x, D') \subseteq A_x$. Then

$$(2.46) \quad \nu_x P^x(X_t = y) = \nu_x P^x(X_t = y, X_{t'} \in A_y) + \nu_x P^x(X_t = y, X_{t'} \in A_x).$$

To bound the first term in (2.46), and writing $\tau = \tau(x, D')$, we have

$$
\begin{aligned}
P^x(X_t &= y, X_{t'} \in A_y) \\
&= P^x(\tau < t', X_{t'} \in A_y, X_t = y) \\
&\le E^x(1_{\{\tau < t'\}} P^{X_\tau}(X_{t-\tau} = y)) \\
&\le P^x(\tau(x, D') < t') \sup_{z \in \partial B(x, D'), u \le t/2} q_{t-u}(z, y) \nu_y \\
&\le c t^{-d/2} P^x(\tau(x, D') < t/2).
\end{aligned}
$$
(2.47)

Similarly, using symmetry, for the second term in (2.46) we have

$$(2.48) \quad \nu_x P^x(X_t = y, X_{t'} \in A_x) = \nu_y P^y(X_t = x, X_{t'} \in A_x) \\ \le c t^{-d/2} P^y(\tau(y, D') < t/2).$$

It remains to verify that we can use Proposition 2.18 to bound the terms $P^z(\tau(z, D') < t/2)$ for $z = x, y$. Writing $c'_i$ for the constants $c_i$ in Proposition 2.18, taking $c_1$ large enough we have $D' \ge c'_1$, and $t' \ge c'_1 D'$. As $(x, V_x)$ is $\lambda$-very good, and $D \ge V_x$, we have that $(z, r)$ is $\lambda$-good for $z \in B(x, 2D)$, and $r \ge (2D)^\eta$. We have $c'_2 t' / D' \ge (2D)^\eta$ provided $t \ge c_8 D^{1+\eta}$, and since $t \ge c_6 D^2 / \log D$ this holds by adjusting the constant $c_1$. So

$$(2.49) \quad P^z(\tau(w, D') < t/2) \le c \exp(-c' D^2 / t) \quad \text{for } z = x, y,$$

and combining the estimates (2.46)–(2.49) completes the proof. □

REMARK 2.20. This theorem does not give any bound for $q_t(x, y)$ when $D < c_1 \vee V_x$ and $t < D^2$. In this case we still have the global upper bound (2.18). In addition the "long-range" bounds in Theorem 2.3, bound $q_t(x, y)$ in terms of $\widetilde{d}(x, y)$, but we do not have a bound in terms of $D$.

The final result of this section is that, under fairly mild additional conditions, functions which are harmonic for the discrete time process $X_n, n \in \mathbb{Z}_+$ are also harmonic for the continuous time process $X_t, t \in \mathbb{R}_+$. At the end of Section 5 we will use this remark to note that the corrector constructed using the discrete time process also gives us a corrector for the continuous time process $X_t$.

Let $\widehat{X}$ be the discrete time process given by $\widehat{X}_n = X_n$, $n \in \mathbb{Z}_+$. Write

$$(2.50) \quad \widehat{\mathcal{L}} f(x) = \sum_y q_1(x, y) \nu_y (f(y) - f(x)).$$

We say $h$ is $\widehat{\mathcal{L}}$ harmonic if the sum in (2.50) converges absolutely for all $x$, and $\widehat{\mathcal{L}}h(x) = 0$ for all $x$. This implies that $(h(\widehat{X}_n), n \in \mathbb{Z}_+)$ is a $P^x$-martingale for each $x \in G$.

For $x \in G$ let $\kappa_x = \mu_x/\nu_x$ be the jump rate out of $x$ by $X$. Set

(2.51) $$A(K) = \{y \in G : \kappa_y \leq K\}.$$

LEMMA 2.21. *Let $\Gamma$ satisfy Assumption 2.6. In addition suppose that there exist $(x_0, r_0)$ such that $(x_0, r_0)$ is $\lambda$-good, and that there exist $R_0$, $K$ such when $R \geq R_0$ then every self avoiding path $\gamma$ from $x_0$ to $B(x_0, R)^c$ contains at least $R^{1/2}$ points in $A(K)$. Let $h : G \to \mathbb{R}$ be $\widehat{\mathcal{L}}$ harmonic, and satisfy the growth condition*

(2.52) $$|h(x)| \leq C_1 + C_1 d(x_0, x)^p$$

*for some $p \in [1, \infty)$. Then $\mathcal{L}h = 0$, so that $h$ is harmonic for $X$.*

PROOF. By Lemma 2.10(b) we have that $(x, r)$ is $2\lambda$-good if $r/2 = r_0 \vee (1 + C_A\lambda)d(x_0, x)$. So by Lemma 2.12 there exists $C_2$ (depending on $C_1$ and $\lambda$) so that if $s \in [0, 1]$,

(2.53) $$E^x|h(X_s)| \leq cr^{d+p} \leq C_2(r_0^{p+d} + d(x_0, x)^p);$$

it follows that $E^x h(X_t)$ is well defined for any $t \geq 0$. Set for $s \in [0, \infty)$

$$h_s(x) = E^x h(X_s).$$

To prove the lemma, it is sufficient to prove that $h = h_s$ for every $s$; this implies that $h(X_t)$ is a continuous time martingale and hence that $\mathcal{L}h = 0$. We have

$$h_{s+1}(x) = E^x(E^{X_s}h(X_1)) = E^x(h(X_s)) = h_s(x),$$

so $s \to h_s$ has period 1. We extend $h_s$ by periodicity to $s \in \mathbb{R}$. Since $E^x h_s(X_1) = E^x h(X_{1+s}) = h_s(x)$, each $h_s$ is $\widehat{\mathcal{L}}$-harmonic. Let

$$k(x) = \sup_{0 \leq s_1 \leq s_2 \leq 1} |h_{s_2}(x) - h_{s_1}(x)|;$$

note that by (2.53) we have

(2.54) $$k(x) \leq 2 \sup_{s \leq 1} E^x|h(X_s)| \leq 2C_2(r^{p+d} + d(x_0, x)^p).$$

Fix $x \in G$, and write $\kappa = \kappa_x$. Write $P_{xy} = \mu_{xy}/\mu_x$ for the jump probabilities of $X$. Then by conditioning on the time of the first jump of $X$, if it occurs in $[0, 1]$, we obtain

$$h_s(x) = e^{-\kappa}h_s(x) + \sum_y P_{xy} \int_0^1 \kappa e^{-\kappa u} h_{s-u}(y) \, du.$$

So

$$h_s(x)(1-e^{-\kappa})$$
$$(2.55) \qquad = \sum_y P_{xy}\left(\kappa \int_0^1 e^{-\kappa}h_{s-u}(y)\,du + \int_0^1 \kappa(e^{-\kappa u}-e^{-\kappa})h_{s-u}(y)\,du\right)$$
$$= \sum_y P_{xy}\left(\kappa \int_0^1 e^{-\kappa}h_u(y)\,du + \int_0^1 \kappa(e^{-\kappa u}-e^{-\kappa})h_{s-u}(y)\,du\right).$$

Then (2.55) implies that

$$(2.56) \qquad k(x)(1-e^{-\kappa_x}) \le \sum_y P_{xy}k(y)(1-(1+\kappa_x)e^{-\kappa_x}).$$

So if $k(x) > 0$, then there exists $y \sim x$ with $k(y) > k(x)$. Further, if $\kappa_x \le K$, then there exists $\delta > 0$ (depending only on $C_M$ and $K$) such that

$$(2.57) \qquad k(y) \ge (1+\delta)k(x) \qquad \text{for some } y \sim x.$$

Suppose now that there exists $x_1$ with $k(x_1) > 0$. Then there exists a non-intersecting infinite path $\gamma_1$ starting at $x_1$ on which $k$ is strictly increasing. Let $\gamma_2$ be a shortest path from $x_0$ to a closest point $y$ on $\gamma$ to $x_0$, and let $D$ be the length of $\gamma_2$. Combining $\gamma_2$ and the infinite segment of $\gamma_1$ starting at $y$, we obtain a path $\gamma = (x_0, z_1, \ldots)$ for which $k(z_n) > 0$ for all $n > D$. Let $R > R_0$, and let $w_R$ be the first point in $\gamma \cap B(x_0, R)^c \cap A(K)$. Then $R_1 = d(x_0, w_R) \ge R$. So, using (2.54), (2.57) and the condition on $A(K)$,

$$2C_2(r_0^{p+d} + R_1^p) \ge k(w_R) \ge (1+\delta)^{R_1^{1/2}-D}k(x_1),$$

which is a contradiction if $R$ is large enough. $\square$

**3. Lower bounds and Harnack inequalities.** Unlike the papers [10, 12, 34] we will need to make explicit use of heat kernel lower bounds in our proof of the invariance principle Theorem 1.1 (see Lemma 5.9).

In this section we specialize to the case when $\Gamma$ is the $d$-dimensional Euclidean lattice, and $\mu_e$ are bond conductances with $\mu_e \ge 1$. We continue to assume that Assumptions 2.1 and 2.6 hold. Note that balls and distance are with respect to the graph distance on $\mathbb{Z}^d$.

We can follow the arguments in Section 5 of [1] fairly closely. First, as $\mu_{xy} \ge 1$ when $x \sim y$, by comparison with the standard Dirichlet form $\mathcal{E}_0$ on $\mathbb{Z}^d$ we have a weighted Poincaré inequality as in [1], Theorem 4.8.

THEOREM 3.1. *Let $B = B(x_0, R)$, $\rho_B(y) = d(y, B^c)$ and $\varphi(x) = R^2(R \wedge \rho_B(y))^2$. Then if $f: B \to \mathbb{R}$,*

$$(3.1) \qquad \inf_a \sum_{x \in B}(f(x)-a)^2\varphi(x)\nu_x \le CR^2 \sum_{x,y \in B}(f(x)-f(y))^2\varphi(x)\wedge\varphi(y)\mu_{xy}.$$

Using this, and the method of Fabes and Stroock [21] we obtain a lower bound of the form $q_t(x,y) \geq ct^{-d/2}$ when $x, y$ are close enough together.

PROPOSITION 3.2. *Let $x_0 \in \mathbb{Z}^d$ and $R \geq c_1$. Then provided*

(3.2) $\quad\quad\quad\quad (z, c_2 R)$ *is $\lambda$-good for all $z \in B(x_0, R)$,*

*we have*

(3.3) $\quad q_t(x_1, x_2) \geq c_1 t^{-d/2} \quad\quad$ *for $x_1, x_2 \in B(x_0, R/2), \frac{1}{8} R^2 \leq t \leq R^2$.*

PROOF. This can be proved using the argument in [1], Proposition 5.1, with only minor changes. Note that we need to show that $\mathbb{P}^{x_1}(X_t \notin B(x_0, 2R/3)) \leq \frac{1}{2}$ when $t = \theta R^2$ and $\theta$ is sufficiently small (see (5.2) and (5.9) in [1]). [There is a missing minus sign in exponential in the last line of (5.2).] This is done using Lemma 2.16, and so to satisfy (2.40) we need (3.2). □

THEOREM 3.3. *Let $x, y \in \mathbb{Z}^d$, $t > 0$, and write $D = d(x,y)$. Suppose that*

(3.4) $\quad\quad\quad\quad\quad\quad t \geq c_1 \vee V_x^2 \vee D^{1+\eta}.$

*Then*

(3.5) $\quad\quad\quad\quad\quad\quad q_t(x, y) \geq c_2 t^{-d/2} e^{-c_3 d(x,y)^2 / t}.$

PROOF. The proof as in [1], Lemma 5.2, Theorem 5.3, follows by a standard chaining argument. We just give the details of the conditions on $V_x$, $D$ and $t$ needed to make this argument work.

First, if $D^2 < t$ then the lower bound in (3.5) is just $t^{-d/2}$, so we can use Proposition 3.2. We set $R = ct^{1/2}$. Then $t \geq V_x^2$ implies $R \geq V_x$, so $B(x, R)$ is $\lambda$-very good, and so as $cR \geq R^\eta$, (3.2) holds.

If $D^2 \geq t$ then we set $R = 2D$, $r = ct/D$. We apply Proposition 3.2 in a chain of balls $B_i = B(z_i, r)$ linking $x$ and $y$. (See [1], Lemma 5.2, or [21] for details of the calculations.) Since $D \geq t^{1/2} \geq V_x$, we have that $(x, R)$ is $\lambda$-very good, and hence that $(z, cr')$ is $\lambda$-good for all $r' \geq R^\eta$, $z \in B(x, R)$. As $r = ct/D \geq c'R^\eta$, (3.2) holds for all the balls $B_i$. □

REMARK 3.4. 1. Note that the lower bounds do not extend to the range when $t \simeq D$. The difficulty is that if $t \simeq D$, then we need Proposition 3.2 for a chain of balls of radius $O(1)$ connecting $x$ and $y$. The hypothesis "very good" is not enough to ensure this.

However, the chaining argument does not need (3.2) for all points in $B(x_0, R)$, but just for a suitable chain connecting $x$ and $y$. In [1] this fact was used to obtain full Gaussian lower bounds. It is likely that the same approach will work for the random conductance model, but we do not pursue this point, since the bounds in Theorem 3.3 are enough for most applications.

2. A well-known theorem (see [23, 33]) states that for Brownian motion on a manifold Gaussian bounds are equivalent to two conditions: volume doubling plus a family of Poincaré inequalities. This theorem was extended to graphs in [18]. Since we have volume doubling (for $\nu$) and the Poincaré inequalities hold (since $\mu_e \geq 1$), one might therefore ask if Theorems 2.19 and 3.3 follow immediately from known results.

However, it is clear that some conditions on $\mu_e$ are needed before Theorem 2.19 holds—one has to prevent $X$ from moving a long distance in a very short time. In fact, examination of the theorems in [18, 23, 33] shows that in each case there is a "hidden" additional assumption which prevents the process from moving too quickly. For example, [18] considers a discrete time nearest neighbour random walk.

For $B \subset \mathbb{Z}^d$ let $q_t^B(x,y)$ be the transition density for the processes $X$ killed on exiting from $B$.

LEMMA 3.5. *Let $(x_0, R)$ be $\lambda$-very good. Then*

(3.6) $\quad q_t^{B(x_0,R)}(x,y) \geq c_1 t^{-d/2}, \qquad x, y \in B(x_0, 3R/4), c_2 R^2 \leq t \leq R^2.$

PROOF. Using Theorem 2.19 and Proposition 3.2, this follows, as in [1], Lemma 5.8, by the argument in [21], Lemma 5.1. □

We now give a parabolic Harnack inequality (PHI) for $X$. The statement requires a little extra notation. If $A \subset \mathbb{Z}^d$ we write $\partial A = \{y : y \sim x$ for some $x \in A\}$ for the exterior boundary of $A$, and $\overline{A} = A \cup \partial A$. We call a function $u(t,x)$ *caloric* in a space–time region $Q = A \times (0, T) \subset [0, \infty) \times \mathbb{Z}^d$ if $u$ is defined on $\overline{Q} = \overline{A} \times [0, T]$ and

$$\frac{\partial}{\partial t} u(t,x) = \mathcal{L}_V u(t,x), \qquad (t,x) \in Q.$$

Write $Q(x,R,T) = B(x,R) \times (0,T]$, $Q_-(x,R,T) = B(x, \tfrac{1}{2}R) \times [\tfrac{1}{4}T, \tfrac{1}{2}T]$ and $Q_+(x,R,T) = B(x, \tfrac{1}{2}R) \times [\tfrac{3}{4}T, T]$.

DEFINITION 3.6. We say the parabolic Harnack inequality (PHI) holds with constant $C_P$ for $Q = Q(x, R, T)$ if whenever $u = u(t, x)$ is nonnegative and caloric on $Q$, then

(3.7) $$\sup_{(t,x) \in Q_-(x,R,T)} u(t,x) \leq C_P \inf_{(t,x) \in Q_+(x,R,T)} u(t,x).$$

THEOREM 3.7 (Parabolic Harnack inequality). *There exists a constant $C_P$ such that if $(x, R)$ is $\lambda$-very good. Then the PHI holds with constant $C_P$ in $Q(x, R, R^2)$.*

PROOF. Using the heat kernel bounds in Theorems 2.19 and 3.5, and Lemma 3.5, this follows by the same argument as in [3], Theorem 3.1. □

**4. Heat kernel bounds for the RCM.** In this section we prove Theorem 1.2. Let $E_d$ be the edges of the Euclidean lattice $\mathbb{Z}^d$, and let $\Omega = [1, \infty]^{E_d}$. Let $\mathbb{P}$ be a probability measure on $\Omega$ which makes the coordinates i.i.d. with a law on $[1, \infty)$. We set $\mu_e(\omega) = \omega(e)$ for $e \in E_d$, and for each $\omega \in \Omega$ we consider the random walk $X$ on the graph $(\mathbb{Z}^d, E_d)$ with conductances $\mu_e(\omega)$.

Using the notation of Section 2 we take $\nu_x = 1$ for all $x$, so that we can take $C_M = 1$. We write $P_\omega^x$ for the law of $X$ started at $x$, and

$$q_t^\omega(x, y) = P_\omega^x(X_t = y)$$

for the transition density of $X$.

LEMMA 4.1. *The graph $(\mathbb{Z}^d, E_d)$, conductances $\mu_e$ and random walk $X$ satisfy Assumptions 2.1 and 2.6 with $\mathbb{P}$-probability 1.*

PROOF. Assumption 2.1 is immediate; note we can take $C_D = 2d$. Setting $C_V = 2^d$ Assumption 2.6(1) is also immediate.

Since we have $\mu(e) \geq 1$ for all edges in $\mathbb{Z}^d$, if $\mathcal{E}_0$ is the Dirichlet form of the standard continuous time random walk on $\mathbb{Z}^d$, then $\mathcal{E}(f, f) \geq \mathcal{E}_0(f, f)$, so that the standard Nash inequality on $\mathbb{Z}^d$ (see [14]) implies Assumption 2.6(2) with a constant $C_N$ depending only on $d$. □

In what follows we set $C_A = 1$; thus none of the constants $C_D, C_A, C_M, C_N, C_V$ depend on the law of $\mu_e$ (apart from the fact that $\mathbb{P}(\mu_e \in [1, \infty)) = 1$).

Let $d(x, y)$ be the graph metric on $(\mathbb{Z}^d, E_d)$, and $\widetilde{d} = \widetilde{d}(\omega)$ be the metric given by (2.6); as in the previous sections we write $\widetilde{B}(x, r)$ for balls in the $\widetilde{d}$ metric. Write $B_E(x, r) = \{y \in \mathbb{R}^d : |x - y| < r\}$ for the Euclidean ball center $x$ and radius $r$.

LEMMA 4.2. *There exists a constant $\lambda_0 > 0$ such that*

(4.1) $$\mathbb{P}(\widetilde{B}(0, r) \subset B_E(0, \lambda_0 r)) \geq 1 - c_1 e^{-c_2 r}.$$

PROOF. We use results on first passage percolation from [25]. As in [25] let $b_{0,n}$ be the first time $\widetilde{B}(0, t)$ reaches the hyperplane $\{x_1 = n\}$. Using [25], Theorem 2.18, there exists $\mu_0$ such that $\lim_n n^{-1} b_{0,n} = \mu_0$, a.s. and in $L^1$. By [25], Theorem 1.15, we have $\mu_0 > 0$. By [25], Theorem 5.2, there exist $c_3, c_4 > 0$ such that

(4.2) $$\mathbb{P}(b_{0,n} < \tfrac{1}{2} n \mu_0) \leq c_3 e^{-c_4 n}, \qquad n \geq 0.$$

The times for $\widetilde{B}(0, t)$ to hit each hyperplane $\{x_i = \pm n\}$, for $i = 1, \ldots, d$ have the same law as $b_{0,n}$, so we deduce

(4.3) $$\mathbb{P}(\widetilde{B}(0, \tfrac{1}{2} \mu_0 n) \subset [-n, n]^d) \geq 1 - 2dc_3 e^{-c_4 n}, \qquad n \geq 0,$$

and (4.1) follows easily. □

Note that $\lambda_0$ does depend on the law of $\mu_e$. We fix $\eta \in (0,1)$, and define *good* and *very good* as in Section 2, with $\lambda$ replaced by $\lambda_0$; and we write $V_x$ for the smallest integer such that $(x, V_x)$ is very good.

THEOREM 4.3. (a)
$$\mathbb{P}((x,r) \text{ is not good}) \leq ce^{-cr}, \qquad r \geq r_0.$$

(b)
$$(4.4) \qquad \mathbb{P}(V_x \geq n) \leq c\exp(-cn^\eta).$$

PROOF. Let $G(y,r) = \{(y,r) \text{ is good}\}$, and $F(R) = \{(y,r) \text{ is good for all } y \in B(0,2R), r \geq R^\eta\}$. Then
$$\mathbb{P}(G(y,r)^c) \leq \sum_{n=r}^\infty ce^{-cn} \leq ce^{-cr}.$$

So,
$$\mathbb{P}(F(R)^c) \leq cR^d \sum_{k=R^\eta}^\infty c_3 e^{-c_2 k} \leq c\exp(-cR^\eta),$$

and since $\{V_0 \geq n\} = \bigcup_n^\infty F(k)^c$, (b) follows. □

Using Lemma 2.11 we obtain the following:

COROLLARY 4.4. *$X$ is conservative with $\mathbb{P}$-probability 1.*

COROLLARY 4.5. *Let $x \in \mathbb{Z}^d$. Then*
$$(4.5) \qquad \lim_{M \to \infty} \lim_{t \to \infty} P_\omega^0(|X_t| \geq M\sqrt{t}) = 0, \qquad \mathbb{P}\text{-}a.s.$$

PROOF. By Lemma 2.15, for $t \geq cV_0^2(\omega)$,
$$P_\omega^0(|X_t| \geq M\sqrt{t}) \leq cM^{-1} t^{-1/2} E_\omega^0 d(0, X_t) \leq cM^{-1},$$

and (4.5) follows. □

THEOREM 4.6. *There exist r.v. $U_x, x \in \mathbb{Z}^d$, such that*
$$(4.6) \qquad \mathbb{P}(U_x(\omega) \geq n) \leq c_1 \exp(-c_2 n^\eta),$$

and if $|x-y| \vee t^{1/2} \geq U_x$, then

(4.7) $\quad q_t^\omega(x,y) \leq c_3 t^{-d/2} e^{-c_4|x-y|^2/t} \qquad$ when $t \geq |x-y|$,

(4.8) $\quad q_t^\omega(x,y) \leq c_3 \exp(-c_4|x-y|(1 \vee \log(|x-y|/t))) \qquad$ when $t \leq |x-y|$.

Further,

(4.9) $\qquad q_t^\omega(x,y) \geq c_6 t^{-d/2} e^{-c_7|x-y|^2/t} \qquad$ if $t \geq U_x^2 \vee |x-y|^{1+\eta}$.

PROOF. We take $U_x = c_8(V_x + 1)$ where $c_8 \geq 1$. The bounds then follow from Theorems 2.19 and 3.3. [Note that the bounds (4.7) and (4.8) are of the same form if $d(x,y) \leq t \leq cd(x,y)$.] We use the constant $c_8$ to adjust between the the Euclidean metric $|x-y|$ and the graph metric $d(x,y)$, and to absorb the conditions $d(x,y) \geq c$ and $t \geq c$ into (4.6). □

THEOREM 4.7. *There exists a constant $C_P$ and r.v. $U_x, x \in \mathbb{Z}^d$ with*

(4.10) $\qquad \mathbb{P}(U_x(\omega) \geq n) \leq c_1 \exp(-c_2 n^\eta),$

*such that if $R \geq U_x$ then a PHI with constant $C_P$ holds for $Q(x, R, R^2)$.*

PROOF. This is immediate from Theorem 3.7 and (4.4). □

The PHI implies an elliptic Harnack inequality (EHI), which holds for the CSRW as well as the VSRW. A function $h$ is harmonic on $A \subset \mathbb{Z}^d$ if it is defined on $\overline{A}$ and $\mathcal{L}_V h(x) = 0$ [or equivalently $\mathcal{L}_C h(x) = 0$] for $x \in A$.

COROLLARY 4.8. *There exists a constant $C_H$ and r.v. $U_x, x \in \mathbb{Z}^d$ with*

(4.11) $\qquad \mathbb{P}(U_x(\omega) \geq n) \leq c_1 \exp(-c_2 n^\eta),$

*such that if $R \geq U_x$, then an EHI with constant $C_H$ holds for $B(x, R)$; if $h \geq 0$ is harmonic in $B(x, R)$, then*

(4.12) $\qquad \sup_{y \in B(x, R/2)} h(y) \leq C_H \inf_{y \in B(x, R/2)} h(y).$

We have the following averaged bounds:

THEOREM 4.9. (a) *Let $x, y \in \mathbb{Z}^d$ and $t \geq c_1 \vee |x-y|^{1+\eta}$. Then*

(4.13) $\qquad c_2 t^{-d/2} e^{-c_3|x-y|^2/t} \leq \mathbb{E} q_t^\omega(x,y) \leq c_4 t^{-d/2} e^{-c_5|x-y|^2/t}.$

(b) *We have*

(4.14) $\qquad \mathbb{E} E_\omega^0 |X_t|^2 \leq c_6 t, \qquad t \geq 1.$

PROOF. (a) Let $D = |x - y|$. Choose $c_1$ so that $\mathbb{P}(U_x > c_1^{1/2}) < \frac{1}{2}$. Then if $t \geq c_1 \vee D^{1+\eta}$, by (4.9),

$$\mathbb{E}q_t^\omega(x,y) \geq \mathbb{E}(q_t^\omega(x,y); U_x^2 < c_1) \geq \tfrac{1}{2}ct^{-d/2}e^{-cD^2/t}.$$

For the upper bound, let $\eta' = 1 - \eta$, and $R'_x$ be the r.v. given in Theorem 4.6 using $\eta'$ instead of $\eta$. Then by (4.7) and (4.6),

$$\mathbb{E}q_t^\omega(x,y) = \mathbb{E}(q_t^\omega(x,y); R'_x > D) + \mathbb{E}(q_t^\omega(x,y); R'_x \leq D)$$

$$\leq ct^{-d/2}e^{-cD^{\eta'}} + ct^{-d/2}e^{-cD^2/t}.$$

Since the second term is larger, the upper bound in (4.13) follows.

(b) We have

$$(4.15) \qquad E_\omega^0|X_t|^2 = \sum_x |x|^2 q_t^\omega(0,x).$$

We split the sum in (4.15) into three parts. First,

$$(4.16) \qquad \sum_{|x| < U_0} |x|^2 q_t^\omega(0,x) \leq U_0^2.$$

Next, using (4.7),

$$(4.17) \qquad \sum_{U_0 < |x| < ct} |x|^2 q_t^\omega(0,x) \leq \sum_{U_0 < |x| < ct} |x|^2 ct^{-d/2} e^{-c|x|^2/t} \leq ct.$$

Finally, using (4.8),

$$(4.18) \qquad \sum_{ct \vee U_0 \leq |x|} |x|^2 q_t^\omega(0,x) \leq \sum_{ct \leq |x|} |x|^2 ce^{-c|x|} \leq c'.$$

Combining (4.15), (4.16) and (4.17) gives

$$E_\omega^0|X_t|^2 \leq ct + c'U_0^2,$$

and as by (4.6) $\mathbb{E}U_0^2 < \infty$ and $t \geq 1$, we obtain (4.14). □

REMARK 4.10. Combining Lemma 2.8, Theorems 4.6 and 4.9 completes the proof of Theorem 1.2.

Now let

$$(4.19) \qquad X_t^{(\varepsilon)} = \varepsilon X_{t/\varepsilon^2}, \qquad 0 < \varepsilon \leq 1.$$

THEOREM 4.11. *Let $T > 0$, $\delta > 0$, $r > 0$. Then*

$$\lim_{R \to \infty} \sup_{\varepsilon} P_\omega^0 \Big(\sup_{s \leq T} |X_s^{(\varepsilon)}| > R\Big) \to 0, \tag{4.20}$$

$$\lim_{\delta \to 0} \limsup_{\varepsilon \to 0} P_\omega^0 \Big(\sup_{|s_1 - s_2| \leq \delta, s_i \leq T} |X_{s_2}^{(\varepsilon)} - X_{s_1}^{(\varepsilon)}| > r\Big) = 0. \tag{4.21}$$

PROOF. By Theorem 4.6, if $R \geq U_0$, then

$$P_\omega^0 \Big(\sup_{s \leq T} X_s \geq R\Big) \leq c \exp(-cR^2/T).$$

So if $R \geq U_0$, then $R/\varepsilon \geq U_0$ and

$$P_\omega^0 \Big(\sup_{s \leq T} |X_s^{(\varepsilon)}| > R\Big) = P_\omega^0 \Big(\sup_{s \leq T/\varepsilon^2} |X_s| \geq R/\varepsilon\Big) \leq c \exp(-cR^2/T),$$

proving (4.20).

To prove (4.21) write

$$p(T, \delta, r) = P_\omega^0 \Big(\sup_{|s_1 - s_2| \leq \delta, s_i \leq T} |X_{s_2} - X_{s_1}| > r\Big), \tag{4.22}$$

so that

$$P_\omega^0 \Big(\sup_{|s_1 - s_2| \leq \delta, s_i \leq T} |X_{s_2}^{(\varepsilon)} - X_{s_1}^{(\varepsilon)}| > r\Big) = p(T/\varepsilon^2, \delta/\varepsilon^2, r/\varepsilon).$$

We begin by bounding $p(T, \delta, r)$ for fixed $T$, $\delta$ and $r$. Let $\kappa \in (0, \frac{1}{2})$, $U_R^* = \max_{x \in B(0,R)} U_x$, and $H(R) = \{U_R^* \leq R^\kappa\}$. Then

$$\mathbb{P}(H(R)^c) \leq cR^d \exp(-cR^{\kappa \eta}), \tag{4.23}$$

so by Borel–Cantelli there exists $R_0 = R_0(\omega)$ such that $\omega \in H(R)$ for all $R \geq R_0$.

Let

$$Z_k = \sup_{0 \leq s \leq \delta} |X_{k\delta + s} - X_{k\delta}|. \tag{4.24}$$

Then if $K = \lfloor T/\delta \rfloor$ and $Z^* = \max_{0 \leq k \leq K} Z_k$, it is enough to control $Z^*$ since

$$\sup_{|s_1 - s_2| \leq \delta, s_i \leq T} |X_{s_2} - X_{s_1}| \leq 2Z^*.$$

Let $R \geq 1$. Then

$$P_\omega^0(Z^* \geq r) \leq P_\omega^0(\tau(0, R) \leq T) + P_\omega^0(Z^* \geq r, \tau(0, R) > T). \tag{4.25}$$

By Proposition 2.18 we have

$$P_\omega^0(\tau(0, R) \leq T) \leq c \exp(-cR^2/T), \tag{4.26}$$

provided that $(0, R)$ is very good. For this it is sufficient that $R \geq U_0(\omega)$. Now,

$$P_\omega^0(Z^* \geq r, \tau(0,R) > T) \leq \sum_{k=0}^{K} P_\omega^0(Z_k \geq r, X_{k\delta} \in B(0,R))$$

$$\leq \sum_{k=0}^{K} \sum_{y \in B(0,R)} P_\omega^y(\tau(y,r) < \delta) P_\omega^0(X_{k\delta} = y).$$

Again by Proposition 2.18, for $y \in B(0, R)$,

(4.27) $$P_\omega^y(\tau(y,r) < \delta) \leq c \exp(-cr^2/\delta),$$

provided $r \geq U_R^*$. This will hold if $R \geq R_0(\omega)$ and $r \geq R^\kappa$. Combining (4.25), (4.26), (4.27), we obtain

(4.28) $$P_\omega^0(Z^* \geq r) \leq c \exp(-cR^2/T) + c(T/\delta) \exp(-cr^2/\delta),$$

provided $R \geq R_0(\omega)$ and $r \geq R^\kappa$.

Hence

(4.29) $$p(T/\varepsilon^2, \delta/\varepsilon^2, 2r/\varepsilon) \leq c \exp(-cR^2/T) + c(T/\delta) \exp(-cr^2/\delta),$$

provided $R > \varepsilon R_0$ and $r \geq R^\kappa \varepsilon^{1-\kappa}$. For fixed $r$, $\delta$ choose $R$ so that $R \geq R_0$ and $R^2/T \geq r^2/\delta$. Then

$$p(T/\varepsilon^2, \delta/\varepsilon, 2r/\varepsilon) \leq cT\delta^{-1} \exp(-cr^2/\delta) \qquad \text{when } \varepsilon^{1-\kappa} \leq rR^{-\kappa}.$$

Hence

$$\limsup_{\varepsilon \to 0} p(T/\varepsilon^2, \delta/\varepsilon, 2r/\varepsilon) \leq cT\delta^{-1} \exp(-cr^2/\delta),$$

and (4.21) follows. □

For $n \in \mathbb{N}$ let $\widehat{X}_n = X_n$, and set

(4.30) $$\widehat{X}_t^{(\varepsilon)} = \varepsilon \widehat{X}_{\lfloor t/\varepsilon^2 \rfloor}, \qquad 0 < \varepsilon \leq 1.$$

LEMMA 4.12. *For any $u > 0$,*

(4.31) $$\lim_{\varepsilon \to 0} P_\omega^0 \left( \sup_{0 \leq s \leq T} |\widehat{X}_s^{(\varepsilon)} - X_s^{(\varepsilon)}| > u \right) = 0.$$

PROOF. In the notation of the previous theorem, it is sufficient to bound $p(T/\varepsilon^2, 1, u/\varepsilon)$; using (4.29) we have

$$p(T/\varepsilon^2, 1, u/\varepsilon) \leq cT\varepsilon^{-2} \exp(-cu^2/\varepsilon^2),$$

provided there exists an $R$ with $R \geq \varepsilon R_0(\omega)$, $R^2 \geq Tu^2 \varepsilon^{-2}$ and $u/\varepsilon \geq R^\kappa$. Setting $R = T^{1/2} u/\varepsilon$, we need $uT^{1/2} \geq \varepsilon^2 R_0(\omega)$, and $u^{1-\kappa} \geq \varepsilon^{1-\kappa} T^{\kappa/2}$, so these bounds hold for all sufficiently small $\varepsilon$. □

**5. Invariance principle.** In this section we prove the invariance principle Theorem 1.1. We assume that the conductances $\mu_e$ are defined on the space $(\Omega, \mathbb{P})$ where

$$\Omega = [1, \infty]^{E_d}.$$

We write $\mu_e(\omega) = \omega(e)$ for the coordinate maps, and make the following assumptions on the environment $(\mu_e)$.

ASSUMPTION 5.1. (1) $(\mu_e)$ is stationary, ergodic, and invariant under symmetries of $\mathbb{Z}^d$.

(2) $\mu_e \in [1, \infty)$ for all $e \in E_d$, $\mathbb{P}$-a.s.

(3) The conclusions of Theorem 1.2 hold for the VSRW associated with $(\mu_e)$.

As explained in the Introduction, our basic approach is to construct the "corrector" $\chi : \Omega \times \mathbb{Z}^d \to \mathbb{R}^d$ so that, for $\mathbb{P}$-a.a. $\omega$ the discrete time process

(5.1) $$\widehat{M}_n = \widehat{X}_n - \chi(\omega, \widehat{X}_n)$$

is a $P_\omega^0$-martingale with respect to the filtration $\widehat{\mathcal{F}}_n = \sigma(\widehat{X}_k, 0 \le k \le n)$.

The key steps in the proof of the invariance principle are:

1. Tightness (a consequence of Theorem 3.5);
2. The invariance principle for the martingale part. This is standard and follows from the ergodicity of our environment (see [29], proof of Theorem 2.6);
3. The almost sure control of the corrector, for which we use the the ergodicity of the environment, the properties (4.11) and (4.12) and the quenched heat kernel estimates in Theorem 1.2 (see [34], or [12], Theorem 2.3). Note that all we need here is the ergodicity of the environment; ergodicity under the action of each direction as stated in [34], Remark 1.3, is not required since one can use the cocycle property of the corrector (see [13, 26]).

We now give the details. Let

$$\Omega_0 = \{\omega : \omega(e) \in [1, \infty) \text{ for all } e\}.$$

Since $\omega(e)$ satisfies Assumption 5.1 we have $\mathbb{P}(\Omega_0) = 1$. We write $\omega = (\omega(e), e \in E_d)$, and $\omega(x, y) = \omega(\{x, y\})$. For $x \in \mathbb{Z}^d$ define $T_x : \Omega \to \Omega$ by

$$T_x(\omega)(z, w) = \omega(z + x, w + x).$$

Let $X$ be the VSRW with generator $\mathcal{L}_V$ given by (1.3), and $q_t^\omega(x, y)$ be the transition density of $X$. As $\nu_x \equiv 1$ is the invariant measure for $X$,

$$q_t^\omega(x, y) = P_\omega^x(X_t = y) = q_t^\omega(y, x).$$

Write

(5.2) $$Q_{xy}(\omega) = q_1^\omega(x,y), \qquad Q_{xy}^{(n)}(\omega) = q_n^\omega(x,y), \qquad x,y \in \mathbb{Z}^d;$$

and note that $Q_{xy} \le 1$ for all $x,y$, with $\sum_y Q_{xy} = 1$. We have

(5.3) $$Q_{xy}^{(n)} \circ T_z = Q_{x+z,y+z}^{(n)}, \qquad Q_{xy}^{(n)} = Q_{yx}^{(n)}.$$

We define the process $Z$, which gives the "environment seen from the particle," by

(5.4) $$Z_t = T_{X_t}\omega, \qquad t \in [0,\infty),$$

and define the discrete time process $\widehat{Z}$ by $\widehat{Z}_n = Z_n$, $n \in \mathbb{Z}_+$.

Let $L^p = L^p(\Omega, \mathbb{P})$. For $F \in L^2$ write $F_x = F \circ T_x$. Then $\widehat{Z}$ has generator

$$\widehat{L}F(\omega) = \sum_{x \in \mathbb{Z}^d} Q_{0x}(\omega)(F_x(\omega) - F(\omega)).$$

Set

$$\widehat{\mathcal{E}}(F,G) = \mathbb{E} \sum_{x \in \mathbb{Z}^d} Q_{0x}(F - F_x)(G - G_x).$$

LEMMA 5.2. *We have for $F \in L^1$,*

(5.5) $$\mathbb{E}F = \mathbb{E}F_x,$$

(5.6) $$\mathbb{E}(Q_{0x}F_x) = \mathbb{E}(Q_{0,-x}F).$$

PROOF. Since $\mathbb{P}$ is invariant by $T_x$ the first relation is immediate. As $(Q_{0,x})_{-x} = Q_{-x,0} = Q_{0,-x}$ by (5.3), $\mathbb{E}(Q_{0x}F_x) = \mathbb{E}((Q_{0x})_{-x}F) = \mathbb{E}(Q_{0,-x}F)$, proving (5.5). $\square$

LEMMA 5.3. *If $F, G \in L^2$, then $\widehat{\mathcal{E}}(F,F) < \infty$, $\widehat{\mathcal{E}}(F,G)$ is defined, and $\widehat{L}F \in L^2$.*

PROOF. Let $F \in L^2$. Then
$$\widehat{\mathcal{E}}(F,F) = \mathbb{E} \sum_{x \in \mathbb{Z}^d} Q_{0x}(F - F_x)^2$$
$$\le 2\mathbb{E} \sum_{x \in \mathbb{Z}^d} Q_{0x}(F^2 + F_x^2)$$
$$= 2\mathbb{E}F^2 + 2\mathbb{E} \sum_{x \in \mathbb{Z}^d} Q_{0x} F_x^2$$
$$= 2\mathbb{E}F^2 + 2\mathbb{E} \sum_{x \in \mathbb{Z}^d} Q_{0,-x} F^2 = 4\|F\|_2^2.$$

Hence $\widehat{\mathcal{E}}(F,G)$ is defined for $F,G \in L^2$. Also, if $F \in L^2$,

$$\mathbb{E}|\widehat{L}F|^2 = \mathbb{E}\sum_{x,y} Q_{0x}Q_{0y}(F_x - F)(F_y - F)$$

$$\leq \mathbb{E}\left[\left(\sum_{x,y} Q_{0x}Q_{0y}(F_x - F)^2\right)^{1/2} \left(\sum_{x,y} Q_{0x}Q_{0y}(F_y - F)^2\right)^{1/2}\right]$$

$$= \widehat{\mathcal{E}}(F,F) \leq 4\|F\|_2^2. \qquad \Box$$

LEMMA 5.4. *Let $F, G \in L^2$. Then*

(5.7) $$\mathbb{E}(G\widehat{L}F) = -\widehat{\mathcal{E}}(F,G).$$

PROOF. Using (5.5) we have

(5.8) $$\mathbb{E}(Q_{0,-x}G(F_{-x} - F)) = \mathbb{E}(Q_{0x}G_x(F - F_x)).$$

So

$$\mathbb{E}(G\widehat{L}F) = \sum_{x \in \mathbb{Z}^d} \mathbb{E}GQ_{0x}(F_x - F)$$

$$= \frac{1}{2}\sum_{x \in \mathbb{Z}^d} \mathbb{E}GQ_{0x}(F_x - F) + \frac{1}{2}\sum_{x \in \mathbb{Z}^d} \mathbb{E}GQ_{0,-x}(F_{-x} - F)$$

$$= \frac{1}{2}\sum_{x \in \mathbb{Z}^d} \mathbb{E}Q_{0x}(GF_x - GF + G_xF - G_xF_x) = -\widehat{\mathcal{E}}(F,G),$$

where we used (5.8) in the last line. $\Box$

Now we look at "vector fields." We define for $G = G(\omega, x): \Omega \times \mathbb{Z}^d \to \mathbb{R}$,

$$\overline{\mathbb{E}}G = \sum_x \mathbb{E}Q_{0x}G(\cdot, x).$$

DEFINITION. We say $G(\omega, x)$ has the *cocycle property* (see [13, 26]) if

(5.9) $$G(T_x\omega, y - x) = G(\omega, y) - G(\omega, x), \qquad \mathbb{P}\text{-a.s.}$$

Let $\mathcal{H} = \overline{L}^2$ be the set of vector fields $G$ with the cocycle property and $\|G\|^2 = \overline{\mathbb{E}}G^2 < \infty$.

LEMMA 5.5. *Let $G = G(\omega, x) \in \overline{L}^2$.*
(a) $G(\omega, 0) = 0$, and $G(T_x\omega, -x) = -G(\omega, x)$.
(b) *If $x_0, x_1, \ldots, x_n \in \mathbb{Z}^d$, then*

(5.10) $$\sum_{i=1}^n G(T_{x_{i-1}}\omega, x_i - x_{i-1}) = G(\omega, x_n) - G(\omega, x_0).$$

PROOF. (a) follows immediately from the definition. For (b), as $G$ has the cocycle property
$$G(T_{x_{i-1}}\omega, x_i - x_{i-1}) = G(\omega, x_i) - G(\omega, x_{i-1}),$$
giving (5.10). □

It is easy to check the following:

LEMMA 5.6.  $\overline{L}^2$ *is a Hilbert space.*

For $F \in L^2$ we set
$$\nabla F(\omega, x) = F(T_x\omega) - F(\omega).$$

LEMMA 5.7.  *If* $F \in L^2$, *then* $\nabla F \in \overline{L}^2$.

PROOF. First,
$$\overline{\mathbb{E}}|\nabla F|^2 = \sum_x \mathbb{E} Q_{0x}(F_x - F)^2 = \widehat{\mathcal{E}}(F, F) < \infty.$$

Also,
$$\nabla F(T_x\omega, y - x) = F(T_{y-x}T_x\omega) - F(T_x\omega)$$
$$= F(T_y\omega) - F(T_x\omega) = \nabla F(\omega, y) - \nabla F(\omega, x),$$
so $\nabla F$ has the cocycle property. □

LEMMA 5.8.  *Let* $G \in \overline{L}^2$. *Then*
$$(5.11) \qquad \mathbb{E} \sum_x Q_{0x}^{(n)} G(\omega, x)^2 \le n\|G\|_2^2.$$

PROOF. Write $a_n^2$ for the left-hand side of (5.11). Then using (5.9),
$$(5.12) \qquad a_n^2 = \mathbb{E} \sum_x \sum_y Q_{0x}^{(n-1)} Q_{xy} (G(T_x\omega, y - x) + G(\omega, x))^2.$$

We now expand the final square in (5.12) and compute the three terms separately. We have
$$\mathbb{E} \sum_x \sum_y Q_{0x}^{(n-1)}(\omega) Q_{xy}(\omega) G(T_x\omega, y - x)^2$$
$$= \mathbb{E} \sum_x \sum_y Q_{-x,0}^{(n-1)}(T_x\omega) Q_{0,y-x}(T_x\omega) G(T_x\omega, y - x)^2$$

$$= \mathbb{E} \sum_x \sum_z Q^{(n-1)}_{-x,0}(\omega) Q_{0,z}(\omega) G(\omega, z)^2$$

$$= \mathbb{E} \sum_z Q_{0,z}(\omega) G(\omega, z)^2 = \|G\|^2.$$

Also,

$$\mathbb{E} \sum_x \sum_y Q^{(n-1)}_{0x}(\omega) Q_{xy}(\omega) G(\omega, x)^2 = \mathbb{E} \sum_x Q^{(n-1)}_{0x}(\omega) G(\omega, x)^2 = a^2_{n-1}.$$

Finally,

$$\mathbb{E} \sum_x \sum_y Q^{(n-1)}_{0x}(\omega) Q_{xy}(\omega) G(\omega, x) G(T_x \omega, y - x)$$

$$= \mathbb{E} \sum_x \sum_z Q^{(n-1)}_{0x}(\omega) Q_{0,z}(T_x \omega) G(\omega, x) G(T_x \omega, z)$$

$$\leq \left( \mathbb{E} \sum_x \sum_z Q^{(n-1)}_{0x}(\omega) Q_{0,z}(T_x \omega) G(\omega, x)^2 \right)^{1/2}$$

$$\times \left( \mathbb{E} \sum_x \sum_z Q^{(n-1)}_{0x}(\omega) Q_{0,z}(T_x \omega) G(T_x \omega, z)^2 \right)^{1/2}$$

$$= a_{n-1} \|G\|.$$

Thus $a_n \leq a_{n-1} + \|G\|$, and so $a_n \leq n\|G\|$. □

Note that the following lemma uses the heat kernel lower bounds.

LEMMA 5.9. *Let $G \in \overline{L}^2$ and $1 \leq p < 2$. Then there exists a constant $c_p < \infty$ such that*

(5.13) $$(\mathbb{E}|G(\cdot, x)|^p)^{1/p} \leq (c_p|x|)\|G\|.$$

*It follows that, $\mathbb{P}$-a.s.,*

(5.14) $$\lim_{n \to \infty} \max_{|x| \leq n} \frac{|G(\omega, x)|}{n^{d+4}} = 0.$$

PROOF. By (5.9) and the triangle inequality we have

$$(\mathbb{E}|G(\cdot, x)|^p)^{1/p} \leq |x|(\mathbb{E}|G(\cdot, e_1)|^p)^{1/p};$$

so it is enough to bound $\mathbb{E}|G(\cdot, e_1)|^p$. By Theorem 1.2 there exists an integer valued random variable $W_0$ with $W_0 \geq 1$ such that $\mathbb{P}(W_0 = n) \leq$

$c_1 \exp(-c_2 n^\delta)$ for some $\delta > 0$ and $q_t^\omega(0, x) \geq c_3 t^{-d/2}$ for $t \geq W_0$. Write $\xi_n = q_n^\omega(0, e_1)$. Then

$$\mathbb{E}|G(\cdot, e_1)|^p = \sum_{n=1}^\infty \mathbb{E}|G(\cdot, e_1)|^p 1_{(W_0 = n)}. \tag{5.15}$$

Let $\alpha = 2/p$, and let $\alpha' = 2/(2-p)$ be its conjugate index. Then using Hölder's inequality and (5.11),

$$\mathbb{E}|G(\cdot, e_1)|^p 1_{(W_0 = n)}$$
$$= \mathbb{E}(\xi_n^{1/\alpha} |G(\cdot, e_1)|^p \xi_n^{-1/\alpha} 1_{(W_0 = n)})$$
$$\leq (\mathbb{E}\xi_n G(\cdot, e_1)^2)^{1/\alpha} (\mathbb{E}\xi_n^{-\alpha'/\alpha} 1_{(W_0 = n)})^{1/\alpha'}$$
$$\leq \left(\mathbb{E}\sum_y Q_{0y}^{(n)} G(0, y)^2\right)^{1/\alpha} ((c_3 n^{-d/2})^{-\alpha'/\alpha} c_1 \exp(-c_2 n^\delta))^{1/\alpha'}$$
$$\leq (n\|G\|^2)^{1/\alpha} c_4 n^{d/2\alpha} \exp(-c_5 n^\delta)$$
$$= c_4 n^{(d+2)/2\alpha} \exp(-c_5 n^\delta) \|G\|^p.$$

Summing the series in $n$ we obtain (5.13).

Using (5.13) with $p = 1$ we have

$$\mathbb{P}\left(\max_{|x| \leq n} |G(\omega, x)| > \lambda_n\right) \leq (2n)^d \max_{|x| \leq n} \mathbb{P}(|G(\omega, x)| > \lambda_n)$$
$$\leq c n^d \lambda_n^{-1} \max_{|x| \leq n} \mathbb{E}|G(\omega, x)| \leq c n^{d+1} \lambda_n^{-1} \|G\|.$$

Taking $\lambda_n = n^{d+3}$ and using Borel–Cantelli gives (5.14). $\square$

Following [29] we introduce an orthogonal decomposition of the space $\overline{L}^2$. Set

$$\overline{L}_p^2 = \mathrm{cl}\{\nabla F, F \in L^2\} \quad \text{in } \mathcal{H},$$

and let $\overline{L}_s^2$ be the orthogonal complement of $\overline{L}_p^2$ in $\mathcal{H}$. (Here $p$ stands for "potential" and $s$ for "solenoidal.")

LEMMA 5.10. *Let $G \in \overline{L}_p^2$. Then for each $x$, $\mathbb{E}G(x, \omega) = 0$.*

PROOF. Fix $x \in \mathbb{Z}^d$. Note first that if $G = \nabla F$, where $F \in L^2$, then $\mathbb{E}G(\omega, x) = \mathbb{E}(F_x - F) = \mathbb{E}F_x - \mathbb{E}F = 0$.

Now let $G \in \overline{L}_p^2$. Then there exist $F_n \in L^2$ such that $G = \lim_n \nabla F_n$ in $\overline{L}^2$. Since $\mathbb{P}(Q_{0x} > 0) = 1$, it follows that $\nabla F_n(\omega, x)$ converges to $G(\omega, x)$ in $\mathbb{P}$-probability. By Lemma 5.9, for each $p \in [1, 2)$ the sequence $\nabla F_n(\omega, x)$ is bounded in $L^p(\Omega, \mathbb{P})$, and therefore $\nabla F_n(\omega, x)$ converges to $G(\omega, x)$ in $L^1(\Omega, \mathbb{P})$. So $\mathbb{E} G(\omega, x) = \lim_n \mathbb{E} \nabla F_n(\omega, x) = 0$. $\square$

We define the semi-direct product measure $\mathbb{P}^* = \mathbb{P} \times P_\omega^0$.

LEMMA 5.11. *Let $G \in \overline{L}_s^2$. Then*

$$\sum_{x \in \mathbb{Z}^d} Q_{0x}(\omega) G(\omega, x) = 0, \qquad \mathbb{P}\text{-}a.s. \tag{5.16}$$

*Hence $M_n = G(\omega, X_n)$ is a $P_\omega^0$-martingale for $\mathbb{P}$-a.a. $\omega$. Further, writing*

$$\|G(\omega, \cdot)\|^2 = \sum_x Q_{0,x}(\omega, x) G(\omega, x)^2,$$

*we have*

$$\langle M \rangle_n = \sum_{k=0}^{n-1} \|G(T_{\widehat{X}_k} \omega, \cdot)\|^2 = \sum_{k=0}^{n-1} \|G(\widehat{Z}_k, \cdot)\|^2. \tag{5.17}$$

*Hence*

$$\mathbb{E}^*(M_n)^2 = n\|G\|^2. \tag{5.18}$$

PROOF. If $F \in L^2$ and $G \in \overline{L}^2$, then using Lemma 5.5,

$$\sum_{x \in \mathbb{Z}^d} \mathbb{E} Q_{0x} G(\omega, x) F_x = \sum_{x \in \mathbb{Z}^d} \mathbb{E} Q_{0x}(T_{-x} \omega) G(T_{-x}\omega, x) F_x(T_{-x}\omega)$$

$$= \sum_{x \in \mathbb{Z}^d} \mathbb{E} Q_{0,-x}(\omega)(-G(\omega, -x)) F(\omega)$$

$$= -\sum_{x \in \mathbb{Z}^d} \mathbb{E} Q_{0x}(\omega) G(\omega, x) F(\omega).$$

Thus

$$\sum_{x \in \mathbb{Z}^d} \mathbb{E} Q_{0x} G(\cdot, x)(F + F_x) = 0. \tag{5.19}$$

If $G \in \overline{L}_s^2$, then

$$0 = \overline{\mathbb{E}}(G \nabla F) = \sum_x \mathbb{E} Q_{0x} G(\cdot, x)(F_x - F),$$

and so $\mathbb{E} \sum Q_{0x} G F = 0$. Since this holds for any $F \in L^2$ we obtain (5.16).

To show that $M$ is a martingale it is enough to prove that for any $x$,

$$E_\omega^0(G(\omega, X_{n+1}) - G(\omega, X_n) | X_n = x) = 0. \tag{5.20}$$

However, using (5.16),

$$E^0_\omega(G(\omega, X_{n+1}) - G(\omega, X_n)|X_n = x)$$
$$= \sum_y Q_{xy}(\omega)(G(\omega, y) - G(\omega, x))$$
$$= \sum_y Q_{0,y-x}(T_x\omega)G(T_x\omega, y - x) = 0.$$

Recall that $\langle M \rangle$ is the unique predictable process so that $M_n^2 - \langle M \rangle_n$ is a martingale. We have

$$E^x_\omega(M^2_{n+1} - M^2_n|\widehat{X}_n = y) = E^x_\omega((M_{n+1} - M_n)^2|\widehat{X}_n = y)$$
$$= \sum_z Q_{yz}(\omega)(G(\omega, z) - G(\omega, y))^2$$
$$= \sum_z Q_{0,z-y}(T_y\omega)(G(z - y, \omega))^2$$
$$= \|G(T_y\omega, \cdot)\|^2,$$

and (5.17) follows.

Finally,

$$\mathbb{E}^* M_n^2 = \mathbb{E}(E^0_\omega M_n^2) = \mathbb{E}(E^0_\omega \langle M \rangle_n) = \sum_{k=0}^{n-1} \mathbb{E}\|G(T_{\widehat{X}_k}\omega, \cdot)\|^2 = n\|G\|^2. \quad \square$$

Let $\Pi : \mathbb{R}^d \to \mathbb{R}^d$ be the identity, and write $\Pi_j$ for the $j$th coordinate of $\Pi$. Then $\Pi_j(y - x) = \Pi_j(y) - \Pi_j(x)$, so $\Pi_j$ has the cocycle property. Further by (4.14),

$$\overline{\mathbb{E}}|\Pi_j|^2 = \mathbb{E}\sum_x Q_{0x}|x_j|^2 < \infty,$$

so $\Pi_j \in \mathcal{H}$. So we can define $\chi_j \in \overline{L}_p^2$ and $\Phi_j \in \overline{L}_s^2$ by

$$\Pi_j = \chi_j + \Phi_j \in \overline{L}_p^2 \oplus \overline{L}_s^2;$$

this gives our definition of the corrector $\chi = (\chi_1, \ldots, \chi_d) : \Omega \times \mathbb{Z}^d \to \mathbb{R}^d$. We will sometimes write $\chi(x)$ for $\chi(\cdot, x)$. Note that conventions about the sign of the corrector differ—compare [34] and [12]. As the environment process is invariant under isometries of $\mathbb{Z}^d$, $\|\Phi_j\| = \|\Phi_1\|$ for each $j = 1, \ldots, d$.

The following proposition summarizes the properties of $\chi$ and $\Phi$.

PROPOSITION 5.12. (a) $\widehat{M}_n = \widehat{X}_n - \chi(\omega, \widehat{X}_n)$ is a $P^0_\omega$-martingale.
(b) For each $x \in \mathbb{Z}^d$, $\chi(\cdot, x) \in L^1$.

(c) *For each $j = 1, \ldots, d$*

$$\mathbb{E} \sum_x Q_{0x}(\omega) |\Phi_j(\omega, x)|^2 = \|\Phi_1\|^2 < \infty.$$

(d) *$\chi$ is sublinear on average; for each $\varepsilon > 0$*

(5.21) $$\lim_{n \to \infty} n^{-d} \sum_{|x| \leq n} 1_{(|\chi(\omega, x)| > \varepsilon n)} = 0, \qquad \mathbb{P}\text{-}a.s.$$

PROOF. (a) and (b) are immediate from Lemmas 5.11 and 5.9, and (c) is immediate from the definition of $\Phi_j$ as a projection in $\overline{L}^2$. Let $e_1$ be the unit vector $e_1 = (1, 0, \ldots, 0)$. By Lemma 5.10 we have $\mathbb{E}\chi(\cdot, e_1) = 0$. So since

(5.22) $$\chi(\omega, n e_1) = \sum_{k=1}^{n} \chi(T_{(k-1)e_1} \omega, e_1)$$

and as $\chi$ has the cocycle property, the ergodic theorem implies that $\lim_n n^{-1} \chi(\omega, n e_1) = 0$ $\mathbb{P}$-a.s., and (d) then follows by the results in Section 6 of [26]. $\square$

LEMMA 5.13. *The processes $Z$ and $\widehat{Z}$ are ergodic under the time shift on the environment space $\Omega$.*

PROOF. This is well known; see [17], Lemma 4.9, and Section 3 of [10] for a careful proof in discrete time. $\square$

PROOF OF THEOREM 1.1. We begin with the VSRW. The arguments are very similar to those in [10, 12, 34], so we only mention the key points. We define

(5.23) $$\widehat{M}_n = \Phi(\omega, \widehat{X}_n), \qquad \widehat{M}_t^{(\varepsilon)} = \varepsilon \widehat{M}_{\lfloor t/\varepsilon^2 \rfloor}, \qquad t \geq 0,$$

so that

(5.24) $$\widehat{X}_t^{(\varepsilon)} = \varepsilon \widehat{X}_{\lfloor t/\varepsilon^2 \rfloor} = \widehat{M}_t^{(\varepsilon)} + \varepsilon \chi(\omega, \varepsilon^{-1} \widehat{X}_t^{(\varepsilon)}).$$

Thus it is sufficient to prove that the martingale $\widehat{M}^{(\varepsilon)}$ converges to a multiple of Brownian motion, and that for $\mathbb{P}$-a.a. $\omega$, the second term in (5.24) converges in $P_\omega^0$-probability to zero.

We start with the control of the corrector, and use [12], Theorem 2.4. This proves that if the corrector $\chi$ has polynomial growth, and is sublinear on average, then Gaussian upper bounds on the heat kernel imply pointwise sublinearity of $\chi$. Thus, using (1.10), (5.14) and (5.21) we have that for $\mathbb{P}$-a.a. $\omega$,

(5.25) $$\lim_{n \to \infty} \max_{|x| \leq n} \frac{|\chi(\omega, x)|}{n} = 0.$$

Given (5.25) the Gaussian upper bounds then imply that, for $\mathbb{P}$-a.a. $\omega$,

(5.26) $\qquad \varepsilon \chi(\omega, \widehat{X}_{\lfloor t/\varepsilon^2 \rfloor}) \to 0 \qquad \text{in } P_\omega^0\text{-probability.}$

For the convergence of $\widehat{M}^{(\varepsilon)}$, we proceed as in [10]. Let $v \in \mathbb{R}^d$ be a unit vector, write $\widehat{M}_n^v = v \cdot M_n$, and let

$$F_K(\omega) = E_\omega^0(|\widehat{M}_1^v|^2; |\widehat{M}_1^v| \geq K).$$

Then $F_K$ is decreasing in $K$ and

$$\mathbb{E} F_K \leq \mathbb{E} F_0 \leq d \|\Phi_1\|^2.$$

In the notation of Lemma 5.11, $F_0(\omega) = \|v \cdot \Phi(\omega, \cdot)\|^2$, and so by (5.17) the covariance process of $\widehat{M}^v$ is

$$\langle \widehat{M}^v \rangle_n = \sum_{k=0}^{n-1} F_0(\widehat{Z}_k).$$

So by Lemma 5.13 we have $n^{-1} \langle \widehat{M}^v \rangle_n \to \mathbb{E} F_0$, $P_\omega^0$-a.s., for $\mathbb{P}$-a.a. $\omega$.

Using the same arguments as in [10], Theorem 6.2, it is straightforward to check the conditions of the Lindeberg–Feller FCLT for martingales (see, e.g., [20], Theorem 7.7.3), and deduce that $v \cdot \widehat{M}^{(\varepsilon)}$ converges to a constant multiple of Brownian motion. Hence $\widehat{M}^{(\varepsilon)}$ converges to an $\mathbb{R}^d$-valued Brownian motion with nonrandom covariance matrix $D$ given by $D_{ij} = \overline{\mathbb{E}} \Phi_i \Phi_j$. Since the law of the random variables $\omega(e)$ is invariant under symmetries of $\mathbb{Z}^d$, we deduce that there exists $\sigma_V^2 \geq 0$ such that $D = \sigma_V^2 I$, and that

(5.27) $\qquad\qquad\qquad \sigma_V^2 = \overline{\mathbb{E}} \Phi_1^2.$

This establishes the convergence of $\widehat{X}^{(\varepsilon)}$; using Lemma 4.12 gives the convergence of $X^{(\varepsilon)}$ to the same limit.

The global upper bounds on $q_t^\omega(0, x)$ in Lemma 2.8 imply that if $\lambda > 0$ and $\lambda t^{1/2} \geq 1$, then

$$P_\omega^0(|X_t| \leq \lambda t^{1/2}) \leq c t^{-d/2} |B(0, \lambda t^{1/2})| \leq c' \lambda^d.$$

Hence there exists $\lambda > 0$ such that for all large $t$,

$$P_\omega^0(|X_t| > \lambda t^{1/2}) \geq \tfrac{1}{2},$$

which implies that $\sigma_V^2 > 0$.

We now consider the CSRW. Recall from (1.1) the definition of $\mu_x(\omega)$, set $F(\omega) = \mu_0(\omega)$, and

(5.28) $\qquad\qquad A_t = \int_0^t \mu_{X_s}\, ds = \int_0^t F(Z_s)\, ds.$

Then if $\tau_t = \inf\{s \geq 0 : A_s \geq t\}$ is the inverse of $A$, the time changed process

(5.29) $$Y_t = X_{\tau_t}$$

is the CSRW.

By the ergodic theorem for the process $Z$,

$$\lim_{t \to \infty} t^{-1} A_t = \mathbb{E} F = 2d \mathbb{E} \mu_e, \qquad \mathbb{P}^*\text{-a.s.}$$

So if $\mathbb{E}\mu_e < \infty$ then $\tau_t/t \to a$ a.s. where $a = 1/2d\mathbb{E}\mu_e > 0$. Let $Y_t^{(\varepsilon)} = \varepsilon Y_{t/\varepsilon^2}$. Then

(5.30) $$Y_t^{(\varepsilon)} = X_{at}^{(\varepsilon)} + (Y_t^{(\varepsilon)} - X_{at}^{(\varepsilon)})$$

and using Theorem 4.11 we have for any fixed $t_0 \geq 0$ that

(5.31) $$\sup_{0 \leq t \leq t_0} |Y_t^{(\varepsilon)} - X_{at}^{(\varepsilon)}|$$

converges in $P_\omega^0$-probability to 0, for $\mathbb{P}$-a.a. $\omega$. Thus $Y^{(\varepsilon)}$ converges to $\sigma_C B_t'$ where $B'$ is a Brownian motion and $\sigma_C^2 = a\sigma_V^2 > 0$.

In the case when $\mathbb{E}\mu_e = \infty$ we have that $\tau_t/t \to 0$, and hence $Y^{(\varepsilon)}$ converges to a degenerate limit. $\square$

We conclude this section by stating a local limit theorem for $q_t(x, y)$ (for the VSRW). Write

$$k_t(x) = (2\pi t \sigma_V^2)^{-d/2} e^{-|x|^2/2\sigma_V^2 t}$$

for the Gaussian heat kernel with diffusion constant $\sigma_V^2$ where $\sigma_V^2$ is as in Theorem 1.1.

THEOREM 5.14. *Let $X$ be the VSRW. Let $T > 0$. For $x \in \mathbb{R}^d$ write $\lfloor x \rfloor = (\lfloor x_1 \rfloor, \ldots, \lfloor x_d \rfloor)$. Then*

(5.32) $$\lim_{n \to \infty} \sup_{x \in \mathbb{R}^d} \sup_{t \geq T} |n^{d/2} q_{nt}^\omega(0, \lfloor n^{1/2} x \rfloor) - k_t(x)| = 0, \qquad \mathbb{P}\text{-a.s.}$$

PROOF. This is proved as in Section 4 of [3]. We have to verify Assumptions 4.1 and 4.4 in [3], but this is straightforward given the invariance principle and heat kernel bounds in Theorems 1.1 and 1.2, and the PHI in Theorem 4.7.

Note that as $\nu$ is the invariant measure for $X$, in Assumption 4.1(d) all we need is that $\nu(\Lambda_n(x,r))/(2n^{1/2}r)d$ converges, and as $\nu_x = 1$ for all $x$; this is easy. (Here $\Lambda_n(x,r) = (xn^{1/2} + [-rn^{1/2}, rn^{1/2}]^d) \cap \mathbb{Z}^d$.) $\square$

REMARK 5.15. In this section we have constructed a corrector $\chi(\omega, x)$ so that the process

$$M_n = X_n - \chi(\omega, X_n), \qquad n \in \mathbb{Z}_+, \tag{5.33}$$

is a (discrete time) martingale. It is natural to ask if $\chi(\omega, \cdot)$ is also a corrector for the continuous time process $X_t$.

For the RCM with i.i.d. conductances it is straightforward to check that the condition in Lemma 2.21 involving the set $A(K)$ holds $\mathbb{P}$-a.s. (see [12], Lemma 3.1, for a similar argument). We can therefore use Lemma 2.21 with $h(\cdot) = \chi(\omega, \cdot)$ to deduce that

$$M_t = X_t - \chi(\omega, X_t), \qquad t \in \mathbb{R}_+, \tag{5.34}$$

is, for $\mathbb{P}$-a.a. $\omega$, a $P_\omega^0$-martingale.

**6. General ergodic environments.** We conclude this paper with some remarks on more general ergodic random environments. First, note that the proof of the invariance principle in Section 5 just uses the facts that the environment is stationary, symmetric and ergodic, and that the heat kernel bounds in Theorem 1.2 hold.

In the proof of Theorem 1.2 the full strength of the assumption that $\mu_e$ were i.i.d. was only used at one point, in Theorem 4.3, where we controlled the probability that a ball was not very good. The heat kernel upper bounds in Section 2 only require Assumptions 2.1 and 2.6, together with a comparison of the metrics $\tilde{d}(x, y)$ and $d(x, y)$. Given these upper bounds, and using the fact that $\mu_e \geq 1$, no additional hypotheses on $\mu_e$ were needed to obtain the lower bounds in Section 3. We therefore have the following:

THEOREM 6.1. *Let $\mu_e$, $e \in E_d$ be a stationary symmetric ergodic environment, satisfying for some $c_1 > 0$,*

$$\mu_e \in [c_1, \infty) \qquad \text{for all } e \in E_d, \mathbb{P}\text{-a.s.} \tag{6.1}$$

*Let $\tilde{d}_\omega(x, y)$ be the metric given by the first passage percolation construction of (2.6), and (as in Definition 2.9) let $V_x(\lambda)$ be the smallest integer such that $(x, V_x(\lambda))$ is $\lambda$-very good. Suppose that there exists $\lambda_0 < \infty$ and $\eta \in (0, 1)$ such that*

$$\mathbb{P}(V_x(\lambda_0) \geq n) \leq c_1 e^{-c_2 n^\eta}. \tag{6.2}$$

*Then the conclusions of Theorems 1.1, 1.2(a)–(c), 1.3, 4.7 and 5.14 all hold for the environment $(\mu_e)$.*

PROOF. We begin by considering the heat kernel bounds in Theorem 1.2. As in Lemma 4.1, it is immediate that Assumptions 2.1 and 2.6 hold for ($\mu_e$), $\mathbb{P}$-a.s. Using the hypothesis (6.2) instead of Theorem 4.3, the arguments in Section 4 [except for Theorem 4.9(a), for which see Remark 6.2 below] hold in this more general context, and give Theorems 1.2 and 4.7.

Given Theorem 1.2, the arguments in Section 5 then give the invariance principle (Theorem 1.1) and local limit theorem (Theorem 5.14).

Combining these results gives the Green function estimates in Theorem 1.3. □

REMARK 6.2. The proof of Theorem 4.9 used the fact that the bounds in Theorem 4.6 hold for $1 - \eta$ as well as for $\eta$. If we only have (6.2) then we obtain

$$(6.3) \qquad \mathbb{E}q_t^\omega(x,y) \leq c_1 t^{-d/2} e^{-c_2|x-y|^2/t} \qquad \text{if } t \geq c_3 \vee |x-y|^{1+\eta},$$

$$(6.4) \qquad \mathbb{E}q_t^\omega(x,y) \geq c_4 t^{-d/2} e^{-c_5|x-y|^2/t} \qquad \text{if } t \geq c_6 \vee |x-y|^{2-\eta}.$$

REMARK 6.3. If $\mu_e$ is bounded and bounded away from 0, so there exist $0 < c_1 \leq c_2 < \infty$ such that $\mathbb{P}(\mu_e \in [c_1, c_2]) = 1$, then the metrics $d(x,y)$ and $\widetilde{d}(x,y)$ are comparable. So, taking $\lambda_0$ large enough, (6.2) holds.

REMARK 6.4. See [12], Lemma 3.1, or [30], Lemma 5.3, for percolation arguments which are more robust than Theorem 4.3 and which may be useful for establishing (6.2) in more general contexts.

REMARK 6.5. If $\mu_e$ is stationary and ergodic, but not invariant with respect to symmetries of $\mathbb{Z}^d$, then if (6.2) holds, we still obtain Theorem 1.2, and the convergence of $X^{(n)}$ to a Brownian motion with covariance matrix $D$. However, $D$ need not be diagonal.

REMARK 6.6. Unlike ergodic bounded conductance models, the results of this paper certainly do not hold for all unbounded stationary symmetric ergodic random environments. For example, let $d = 2, 3, 4$ and let $\mathcal{T}$ be a uniform spanning tree on $\mathbb{Z}^d$ (see [9]). Then $\mathcal{T}$ is 1-sided, so from each $x \in \mathbb{Z}^d$ there is a unique self-avoiding path $\gamma_x$ to infinity. Let $a(x)$ be the first point on this path. Then $a: \mathbb{Z}^d \to \mathbb{Z}^d$ and the path $\gamma_x$ is $\{x, a(x), a^2(x), \ldots\}$.

Let $N(x)$ be the set of points in $\mathcal{T}$ which are disconnected from infinity by deleting the bond $\{x, a(x)\}$, and let $n(x) = |N(x)|$. As $x \in N(x)$, $n(x) \geq 1$ for all $x$. Let $\mu_e = 1$ for edges $e \in E_d$ which are not in $\mathcal{T}$. Each edge $e \in \mathcal{T}$ is of the form $e = \{x, a(x)\}$ for some $x$, set

$$\mu_{\{x,a(x)\}} = n(x) e^{n(x)^2}.$$

Let $T_i$, $i \geq 1$, be the jump times of the VSRW $X$. Then

$$(6.5) \qquad P_\omega^x(X_{T_1} \neq a(x)) = \frac{\sum_{y \neq a(x)} \mu_{xy}}{\mu_{x,a(x)} + \sum_{y \neq a(x)} \mu_{xy}}.$$

Fix $x$, and let the neighbors of $x$ in $\mathcal{T}$ be $a(x), y_1, \ldots, y_k$. Then

$$\sum_{y \neq a(x)} \mu_{xy} = (2d - k - 1) + \sum_{i=1}^{k} n(y_i) e^{n(y_i)^2}.$$

Since $\mu_{x,a(x)} = n(x) e^{n(x)^2}$ and $n(x) = 1 + \sum n(y_i)$, it is easy to see that

$$p_0(x) = P_\omega^x(X_{T_1} \neq a(x)) \leq 2d e^{-n(x)^2} + \max_i e^{n(y_i)^2 - n(x)^2} \leq 2d e^{-n(x)^2} + e^{-n(x)/d}.$$

So $\sum_k p_0(a^k(x)) < \infty$, and it follows that ultimately the process $X$ moves to infinity along a path $\gamma_x$ for some $x$. Since $\sum_k \mu^{-1}_{a^k(x), a^{k+1}(x)} < \infty$ this takes finite time. Hence the quenched invariance principle Theorem 1.1 fails, as well as the Gaussian bounds in Theorem 1.2.

**Acknowledgments.** The authors thank J. Černý and T. Kumagai for valuable discussions.

DEPARTMENT OF MATHEMATICS
UNIVERSITY OF BRITISH COLUMBIA
VANCOUVER, BC V6T 1Z2
CANADA
E-MAIL: barlow@math.ubc.ca

FACHBEREICH MATHEMATIK
TECHNISCHE UNIVERSITÄT BERLIN
STRASSE DES 17. JUNI 136
D-10623 BERLIN
GERMANY
E-MAIL: deuschel@math.tu-berlin.de